\documentclass[12pt]{amsart}

\usepackage[dvips]{graphicx}
\usepackage{amsfonts,amssymb,amscd,bm}
\newcommand{\R}{{\mathbb R}}

\newcommand{\Ae}{\mathcal{A}_e}
\newcommand{\Az}{\mathcal{A}_z}

\newcommand{\Z}{{\mathbb Z}}

\newcommand{\T}{{\mathbb T}}

\newcommand{\e}{{\varepsilon}}

\newcommand{\be}{{\beta}}

\newcommand{\La}{{\Lambda}}
\newcommand{\ga}{{\gamma}}
\newcommand{\Ga}{{\Gamma}}

\newcommand{\de}{{\delta}}

\newcommand{\De}{{\Delta}}
\newcommand{\si}{{\sigma}}
\newcommand{\Si}{{\Sigma}}
\newcommand{\Om}{{\Omega}}

\newcommand{\ti}{\tilde}

\newcommand{\bu}{\bullet}

\newcommand{\odd}{{\mathrm{Odd}\,}}
\newcommand{\even}{{\mathrm{Even}\,}}

\newcommand{\pr}{{\mathrm{pr}}}

\newcommand{\tg}{{\mathrm{TG}}}

\newcommand{\xy}{{\mathrm{XY}}}

\newcommand{\xa}{{\mathrm{X}}}

\newcommand{\const}{{\mathrm{const}}}

\newcommand{\bd}{{\partial}}

\newcommand{\codim}{\mathrm{codim}\,}

\newcommand{\lra}{\leftrightarrow}

\newcommand{\ed}{{\hfill $\blacksquare$}}

\newcommand{\sk}{\mathrm{SL}}
\newcommand{\SL}{\mathrm{SL}}

\newcommand{\cl}{\mathrm{CL}}

\newcommand{\ds}{\mathrm{DS}}

\newcommand{\gd}{\mathrm{GD}}

\newcommand{\rot}{\mathrm{rot}}

\newcommand{\dxy}{\mathrm{D}_{xy}}

\newcommand{\rxz}{\mathrm{A}_{xz}}
\newcommand{\axz}{\mathrm{A}_{xz}}

\newcommand{\deter}[4]{
 \left| \begin{array}{cc}
  #1 & #3 \\
  #2 & #4
 \end{array} \right| }

\newcommand{\tdotx}{ \begin{picture}(7,8)(0,0)
  \put(1,8){\circle*{1}}
  \put(3,8){\circle*{1}}
  \put(5,8){\circle*{1}}
  \put(0,0){$x$}
 \end{picture} }
\newcommand{\tdoty}{ \begin{picture}(7,8)(0,0)
  \put(1,8){\circle*{1}}
  \put(3,8){\circle*{1}}
  \put(5,8){\circle*{1}}
  \put(0,0){$y$}
 \end{picture} }
\newcommand{\tdotz}{ \begin{picture}(7,8)(0,0)
  \put(1,8){\circle*{1}}
  \put(3,8){\circle*{1}}
  \put(5,8){\circle*{1}}
  \put(0,0){$z$}
 \end{picture} }

\newcommand{\crit}{ \begin{picture}(10,10)(0,0)
  \put(5,0){\oval(8,8)[t]}
  \put(5,0){\line(0,1){8}}
 \end{picture} }
\newcommand{\critver}{ \begin{picture}(10,8)
  \put(0,4){\line(1,0){4}}
  \put(5,4){\circle{3}}
  \put(7,4){\line(1,0){3}}
 \end{picture} }

\newcommand{\cube}{ \begin{picture}(10,10)(0,2)
  \put(5,1){\oval(8,8)[lt]}
  \put(5,9){\oval(8,8)[rb]}
 \end{picture} }
\newcommand{\cubicsn}{ \begin{picture}(8,10)(0,2)
  \put(3,0){\line(0,1){10}}
  \put(-1,5){\oval(8,8)[br]}
  \put(7,5){\oval(8,8)[tl]}
 \end{picture} }
\newcommand{\cusp}{ \begin{picture}(10,10)(0,0)
  \put(1,0){\oval(8,12)[tr]}
  \put(9,0){\oval(8,12)[tl]}
 \end{picture} }
\newcommand{\cuspint}{ \begin{picture}(12,10)(0,2)
  \put(0,1){\line(2,1){10}}
  \put(1,4){\oval(8,10)[tr]}
  \put(9,4){\oval(8,10)[tl]}
 \end{picture} }
\newcommand{\degcusp}{ \begin{picture}(12,10)(0,2)
  \qbezier(2,10)(10,10)(10,0)
  \qbezier(2,6)(8,6)(10,0)
 \end{picture} }

\newcommand{\hangver}{ \begin{picture}(10,8)
  \put(0,4){\line(1,0){8}}
  \put(1,4){\circle*{3}}
 \end{picture} }

\newcommand{\horcusp}{ \begin{picture}(8,10)(0,2)
  \put(0,9){\oval(14,8)[br]}
  \put(0,1){\oval(14,8)[tr]}
 \end{picture} }
\newcommand{\hortang}{ \begin{picture}(10,10)(0,2)
  \put(5,1){\oval(8,8)[t]}
  \put(5,9){\oval(8,8)[b]}
 \end{picture} }
\newcommand{\hortrip}{ \begin{picture}(12,10)(0,2)
  \put(5,1){\oval(8,8)[t]}
  \qbezier(3,1)(5,5)(7,9)
  \qbezier(7,1)(5,5)(3,9)
 \end{picture} }
\newcommand{\inctrip}{ \begin{picture}(16,12)(0,2)
  \qbezier(-1,0)(7,14)(15,0)
  \qbezier(5,1)(7,7)(9,13)
  \qbezier(2,1)(7,6)(12,12)
 \end{picture} }

\newcommand{\maximax}{ \begin{picture}(15,10)(0,2)
  \qbezier(0,2)(3,12)(6,2)
  \qbezier(8,2)(11,12)(14,2)
 \end{picture} }
\newcommand{\maximin}{ \begin{picture}(12,10)(0,0)
  \qbezier(0,0)(3,10)(6,0)
  \qbezier(5,10)(8,0)(11,10)
 \end{picture} }

\newcommand{\minimin}{ \begin{picture}(14,10)(0,2)
  \qbezier(0,8)(3,-4)(6,8)
  \qbezier(8,8)(11,-4)(14,8)
 \end{picture} }

\newcommand{\maxtang}{ \begin{picture}(15,10)(0,0)
  \qbezier(0,0)(6,12)(12,0)
  \qbezier(3,0)(6,12)(9,0)
 \end{picture} }
\newcommand{\quadrup}{ \begin{picture}(12,10)(0,2)
  \put(3,0){\line(2,5){4}}
  \put(0,3){\line(5,2){10}}
  \put(7,0){\line(-2,5){4}}
  \put(0,7){\line(5,-2){10}}
 \end{picture} }

\newcommand{\revcrit}{ \begin{picture}(10,10)(0,0)
  \put(5,8){\oval(8,8)[b]}
  \put(5,0){\line(0,1){8}}
 \end{picture} }

\newcommand{\tang}{ \begin{picture}(10,10)(0,2)
  \put(1,5){\oval(8,8)[r]}
  \put(9,5){\oval(8,8)[l]}
 \end{picture} }
\newcommand{\tangint}{ \begin{picture}(12,10)(0,2)
  \put(0,3){\line(2,1){10}}
  \put(1,5){\oval(8,8)[r]}
  \put(9,5){\oval(8,8)[l]}
 \end{picture} }
\newcommand{\tangver}{ \begin{picture}(8,8)
  \put(6,4){\oval(8,8)[l]}
  \put(2,4){\circle*{2}}
 \end{picture} }
\newcommand{\tantrip}{ \begin{picture}(10,10)(0,2)
  \put(1,10){\line(1,0){8}}
  \put(5,0){\line(0,1){8}}
  \qbezier(1,1)(5,4)(9,8)
  \qbezier(9,1)(5,4)(1,8)
 \end{picture} }
\newcommand{\trip}{ \begin{picture}(11,10)(0,2)
  \put(0,0){\line(1,1){10}}
  \put(5,0){\line(0,1){10}}
  \put(0,10){\line(1,-1){10}}
 \end{picture} }

\newcommand{\tripver}{ \begin{picture}(11,10)(0,2)
  \put(0,0){\line(1,1){10}}
  \put(0,5){\line(1,0){10}}
  \put(0,10){\line(1,-1){10}}
  \put(5,5){\circle*{3}}
 \end{picture} }

\newtheorem{theorem}{Theorem}[section]

\newtheorem{proposition}[theorem]{Proposition}
\newtheorem{lemma}[theorem]{Lemma}

\theoremstyle{definition}
\newtheorem{definition}[theorem]{Definition}
\newtheorem{example}[theorem]{Example}
\newenvironment{sketch}{\noindent\emph{Sketch.}}

\title{A 1-parameter approach to links in a solid torus}

\author[fiedler]{T.~Fiedler}
\address{ Laboratoire Emile Picard, Universit\'e Paul Sabatier,
118 route Narbonne, 31062 Toulouse, France}
\email{ fiedler@picard.ups-tlse.fr }

\author[kurlin]{V.~Kurlin}
\address{ Department of Mathematical Sciences,
Durham University,
Durham DH1 3LE, United Kingdom}
\email{ vitaliy.kurlin@durham.ac.uk }

\subjclass[2000]{57R45, 57M25, 53A04}
\keywords{Knot, braid, singularity, bifurcation diagram, trace graph, 
 diagram surface, canonical loop, trihedral move, tetrahedral move}
\date{ October 18, 2008, the last version
 is available on www.durham.ac.uk/$\sim$dma0vk}

\begin{document}

\begin{abstract}
To an oriented link in a solid torus we associate a trace graph
 in a thickened torus in such a way that links
 are isotopic if and only if their trace graphs
 can be related by moves of finitely many standard types.
The key ingredient is a study of codimension~2 singularities of link diagrams.
For closed braids with a fixed number of strands,
 trace graphs can be recognized up to equivalence 
 excluding one type of moves in polynomial time
 with respect to the braid length.
\end{abstract}

\maketitle



\section{Introduction}
\label{sect:Introduction}

\subsection{Motivation and summary}
\label{subs:Summary}
\noindent
\smallskip

The classical Reidemeister theorem says that plane diagrams represent isotopic links
 in 3-space if and only if they can be related by finitely many moves of 3 types
 corresponding to the codimension~1 singularities of links diagrams, namely
 a triple point $\trip$, simple tangency $\tang$ and ordinary cusp $\cusp$.

We establish the higher order Reidemeister theorem considering
 a canonical 1-parameter family of links in a solid torus and studying
 codimension~2 singularities of resulting link diagrams. 
The 1-parameter family of link diagrams is encoded by a new combinatorial object,
 a trace graph in a thickened torus in such a way that
\emph{trace graphs determine families of isotopic links
 if and only if they can be related by a finite sequence 
 of the 11 moves in Figure~\ref{fig:MovesTraceGraphs} }, 
 see Theorem~\ref{thm:MovesTraceGraphs}.
\smallskip

The conjugacy problem for braids is equivalent to 
 the isotopy classification of closed braids in a solid torus.
Braids are \emph{conjugate} if and only if
 the trace graphs of their closures are \emph{equivalent} 
 through only tetrahedral moves and trihedral moves in 
 Figure~\ref{fig:MovesTraceGraphs}i, \ref{fig:MovesTraceGraphs}xi.
Trace graphs of closed braids can be recognized up to isotopy in a thickened torus
 and trihedral moves in polynomial time with respect to the braid length, 
 see Theorem~\ref{thm:RecognizeUpToThihedralMoves}.
The method provides a new geometric approach to the conjugacy problem 
 for braid groups $B_n$, which still has no efficient solution for $n\geq 5$ strands,
 ie with a polynomial complexity in the braid length.
Very promising steps towards a polynomial solution
 were made by Birman, Gebhardt, Gonz\'alez-Meneses \cite{BGG} and Ko, Lee \cite{KL}.
A clear obstruction is that the number of different conjugacy classes
 of braids grows exponentially even in $B_3$, see Murasugi \cite{Mur}.
\smallskip

Usually links are studied in terms of braids using the theorems of 
 Alexander and Markov, see Birman \cite{Bir}.
The 1-parameter approach is a geometric alternative to the algebraic one:
 conjugacy of braids and Markov moves are replaced
 by a stronger notion of link isotopy and extreme tangency moves in 
 Figure~\ref{fig:MovesTraceGraphs}viii, respectively.
\smallskip


\subsection{Basic definitions}
\label{subs:BasicDefinitions}
\noindent
\smallskip

\begin{figure}[!h]
\includegraphics[scale=1.0]{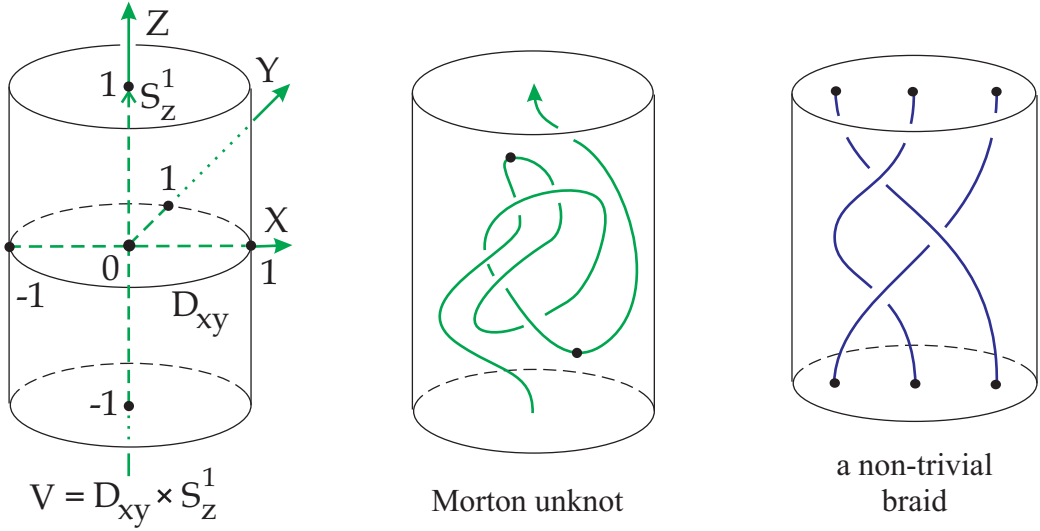}
\caption{Basic notations and examples}
\label{fig:NotationsExamples}
\end{figure}

We work in the $C^{\infty}$-smooth category.
Fix Euclidean coordinates $x,y,z$ in $\R^3$.
Denote by $\dxy$ the unit disk with centre at the origin
 of the horizontal plane $\xy$.
Introduce the \emph{solid torus}
 $V=\dxy\times S_z^1$, where the oriented circle $S_z^1$
 is the segment $[-1,1]_z$ with the identified endpoints, 
 see the left picture of Figure~\ref{fig:NotationsExamples}.
\smallskip

\begin{definition}
\label{def:LinkEquivalence}
An \emph{embedding} is a diffeomorphism onto its image.
An \emph{oriented link} $K\subset V$ is the image
 of an embedding $f:\sqcup_{j=1}^m S_j^1\to V$.
An \emph{isotopy} between two oriented links
 $K_0$ and $K_1$ in $V$ is a smooth map
 $F:(\sqcup_{j=1}^m S_j^1)\times[0,1]\to V$ such that
 $f_0(\sqcup_{j=1}^m S_j^1)=K_0$, $f_1(\sqcup_{j=1}^m S_j^1)=K_1$ and
 the maps $f_r=F(*,r):\sqcup_{j=1}^m S_j^1\to V$ 
 are smooth embeddings for all $r\in[0,1]$.
\medskip

\noindent
Mark $n$ points $p_1,\dots,p_n\in\dxy$.
A \emph{braid} $\be$ on $n$ \emph{strands} is
 the image of a smooth embedding of
 $n$ segments into $\dxy\times[-1,1]_z$
 such that (see Figure~\ref{fig:NotationsExamples})
\smallskip

\noindent
$\bu$
 the strands of $\be$ are monotonic
 with respect to $\pr_z:\be\to S_z^1$;
\smallskip

\noindent
$\bu$
 the lower and upper endpoints of $\be$ are
 $\cup(p_i\times\{-1\})$, $\cup(p_i\times\{1\})$, respectively.
\medskip

\noindent
Braids are considered up to isotopy in the cylinder $\dxy\times[-1,1]_z$, 
 fixed on its boundary.
The isotopy classes of braids form the group denoted by $B_n$.
The \emph{trivial} braid consists of $n$ vertical
 segments $\sqcup_{i=1}^n(p_i\times[-1,1]_z)$.
A braid $\be\in B_n$ is \emph{pure} if the induced permutation $\tilde\be\in S_n$
 is its endpoints is trivial.
The \emph{closed} braid $\hat\be\subset V$ is obtained from 
 $\be\subset\dxy\times[-1,1]_z$ by identifying the bases $\{z=\pm 1\}$.
\ed
\end{definition}
\smallskip

The smoothness of a link $K$ implies
 that the tangent vector $\dot f(s)$ never vanishes on $K$.
The \emph{standard} unknot is given by the trivial embedding
 $S_z^1\to V=\dxy\times S_z^1$.
We introduce a new equivalence relation,
 \emph{strong isotopy}, for links in a solid torus.
For closed braids, the usual isotopy through closed braids is strong.
\smallskip

\begin{definition}
\label{def:StrongIsotopy}
An \emph{extreme pair} of a link $K\subset V$ is
 a pair of either 2 local maxima $\maximax$
 or 2 local minima $\minimin$ of the projection
 $\pr_z:K\to S_z^1$ with the same $z$-coordinate.
A smooth isotopy $F:(\sqcup_{j=1}^m S_j^1)\times[0,1]\to V$
 of links is called \emph{strong} if
 all links in the family $K_r=F(\sqcup_{j=1}^m S_j^1,r)\subset V$
 have {\bf no} extreme pairs for $r\in[0,1]$.
\ed
\end{definition}
\smallskip

H.~Morton proposed the trivial knot in the middle picture of 
 Figure~\ref{fig:NotationsExamples}.
The \emph{Morton} unknot
 is not strongly isotopic to the \emph{standard} unknot $S_z^1$.
The arc between the marked extrema is a long trefoil that
 can not be unknotted by strong isotopy since 
 the marked extrema remain the highest and lowest critical points.
\smallskip


\subsection{Trace graphs of links}
\label{subs:TraceGraphLink}
\noindent
\smallskip

Links are usually represented by plane diagrams 
 with double crossings.
A classical approach to the classification of links
 is to use isotopy invariants, ie functions defined
 on plane diagrams and invariant under the Reidemeister moves.
The Reidemeister moves in Figure~\ref{fig:ReidemeisterMoves} 
 correspond to simplest singularities 
 that can appear in diagrams of links under isotopy, eg
 Reidemeister move~III describes the change of a diagram 
 when a transversal triple intersection $\trip$ appears in the projection.
\smallskip

\begin{figure}[!h]
\includegraphics[scale=1.0]{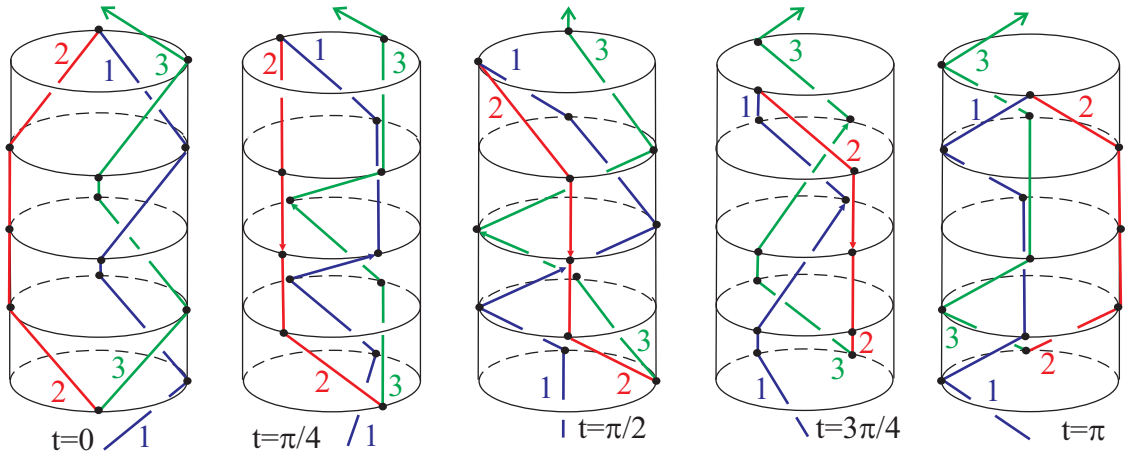}
\caption{Diagrams of rotated trefoils $\rot_t(K)\subset V$ for $t\in[0,\pi]$}
\label{fig:RotatedTrefoils}
\end{figure}

\begin{figure}[!h]
\includegraphics[scale=1.0]{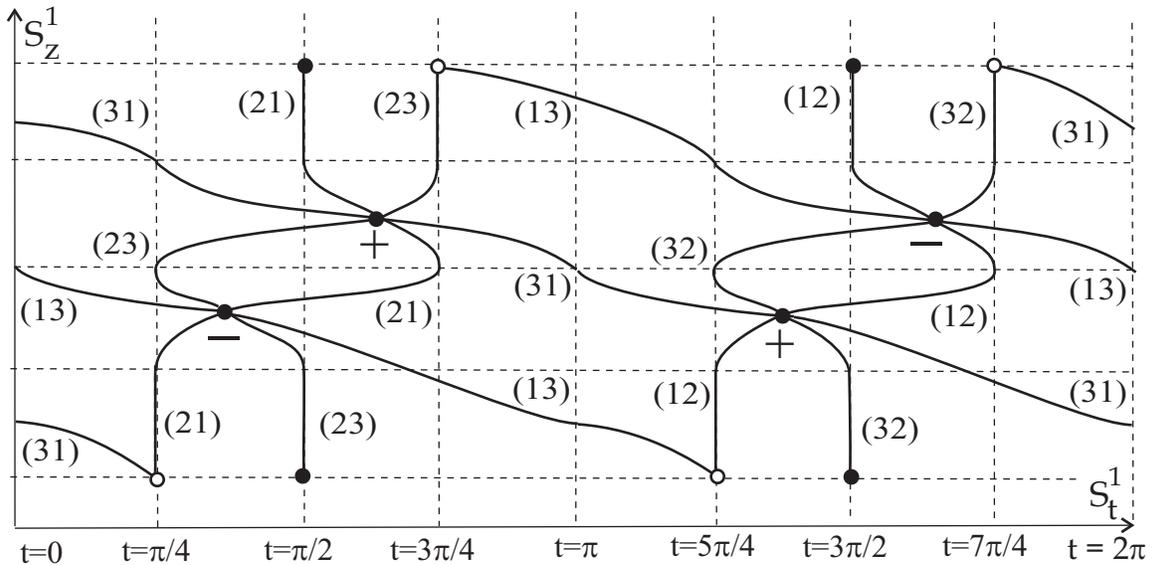}
\caption{The trace graph $\tg(K)$ obtained from the diagrams in Figure~\ref{fig:RotatedTrefoils}.}
\label{fig:TraceGraphTrefoil}
\end{figure}

The analogue of a plane diagram in the 1-parameter approach
 is a 1-parameter family of diagrams 
 obtained by rotating a link in $V$ around $S_z^1$.
This is a 2-dimensional surface
 containing more information about the link
 than only one plane diagram.
The link will be reconstructed up to smooth equivalence 
 from the self-intersection of the surface, the \emph{trace graph}.
Denote by $\rot_t:V\to V$ the \emph{rotation} of the torus $V$ around $S_z^1$ 
 through an angle $t\in[0,2\pi)$, see Figure~\ref{fig:RotatedTrefoils}.
Here $t$ is the parameter on
 \emph{the time circle} $S_t^1=\R/2\pi\Z$ of length $2\pi$.
Let $\axz$ be the vertical annulus $[-1,1]_x\times S_z^1$
 in the solid torus $V$.
Define the \emph{thickened torus} $\T=\axz\times S_t^1$.
We illustrate the rotation of $V$ using the piecewise linear trefoil $K\subset V$
 in Figure~\ref{fig:RotatedTrefoils}, which can be easily smoothed.
Diagrams of rotated trefoils $\rot_t(K)\subset V$ under the orthogonal \emph{projection} 
 $\pr_{xz}:V\to\axz$ are shown in Figure~\ref{fig:RotatedTrefoils}.

\begin{definition}
The \emph{trace graph} $\tg(K)\subset\T$
 of a link $K\subset V$ consists of
 the crossings of the diagrams
 $\pr_{xz}(\rot_t(K))$ over all $t\in S_t^1$.
Mark the intersection points from $K\cap(\dxy\times\{\pm 1\})\subset V$ 
 and also mark each local extremum of $K$ with respect to $\pr_z:K\to S_z^1$.
If $K$ has $m$ components, in general position, the $i$-th component 
 decomposes into $n_i$ vertically monotonic arcs labelled by $A_{iq}$, $q=1,\dots,n_i$.
The 3 monotonic arcs in Figure~\ref{fig:RotatedTrefoils} are numbered simply by 1, 2, 3.
\smallskip

\noindent
Take a point $p\in\tg(K)$, which is a crossing
 of $A_{iq}$ over $A_{js}$ in the diagram
 $\pr_{xz}(\rot_t(K))$ for some $t\in S_t^1$.
Associate to $p$ the ordered \emph{label} $(q_i s_j)$.
Then the edges of the graph $\tg(K)$ are labelled.
In the case $m=1$ we miss the indices $i,j$ of components of $K$
 and label edges by ordered pairs $(qs)$, see Figure~\ref{fig:TraceGraphTrefoil}.
\ed
\end{definition}

For each $t\in S_t^1$, watch the crossings of 
 the diagram $\pr_{xz}(\rot_t(K))\subset\axz\times\{t\}$,
 eg the initial diagram $\pr_{xz}(K)\subset\axz$ at $t=0$ has 
 3 double crossings, which evolve during the rotation of $K$.
At $t=\pi/4$, the lowest crossing becomes a critical crossing $\revcrit$ 
 corresponding to a critical vertex $\critver$ of $\tg(K)$.
At the same $t=\pi/4$ a couple of crossings is born after Reidemeister move II
 associated to a tangency $\tang$.
At $t=\pi/2$ a new crossing is born from a cusp $\cusp$
 after Reidemeister move I, which leads to a hanging vertex $\hangver$ of $\tg(K)$. 
The 2 triple vertices of $\tg(K)$ for $t\in(0,\pi)$ correspond to 
 2 Reidemeister moves III happening during the rotation of $K$.
A combinatorial explicit construction of the trace graph is in 
 Lemma~\ref{lem:Construction}.

\begin{theorem}
\label{thm:MovesTraceGraphs}
Links $K_0,K_1\subset V$ are isotopic in the solid torus $V$ if and only if 
 their labelled trace graphs $\tg(K_0),\tg(K_1)\subset\T$
 can be obtained from each other by an isotopy in $\T$
 and a finite sequence of moves in Figure~\ref{fig:MovesTraceGraphs}.
\end{theorem}
\smallskip

Trace graphs of closed braids have combinatorial features,
 allowing us to recognize them up to all but one type of moves.
The following result of \cite{FK} is one of very few known polynomial 
 algorithms recognizing topological objects up to isotopy.

\begin{theorem}
\label{thm:RecognizeUpToThihedralMoves}
Let $\be,\be'\in B_n$ be braids of length $\leq l$.
There is an algorithm of complexity $C(n/2)^{n^2/8}(6l)^{n^2-n+1}$
 to decide whether $\tg(\hat\be)$ and $\tg(\hat\be')$
 are related by isotopy in $\T$ and trihedral moves, 
 the constant $C$ does not depend on $l$ and $n$.
In the case of pure braids, the power $n^2/8$ can be replaced by 1.
If the closure of a braid is a single curve in the solid torus,
 then the complexity reduces to $Cn(6l)^{n-1}$.
\end{theorem}


\subsection{Scheme of proofs}
\label{subs:SchemeProofs}
\noindent
\smallskip

\begin{figure}[!h]
\includegraphics[scale=1.0]{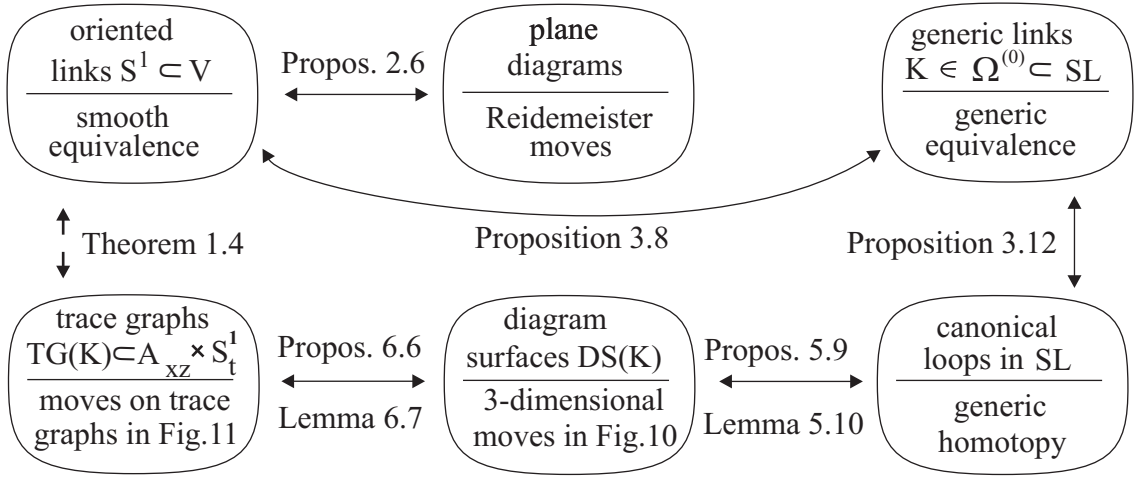}
\caption{A scheme to prove Theorem~\ref{thm:MovesTraceGraphs}.}
\label{fig:SchemeProofs}
\end{figure}

\noindent
The \emph{first} double arrow in Figure~\ref{fig:SchemeProofs} is
 a classical reduction of an equivalence of links to 
 extended Reidemeister moves on plane diagrams in Figure~\ref{fig:ReidemeisterMoves}.
\smallskip

\noindent
The \emph{second} arrow is a new reduction to generic links
 and generic equivalence defined in terms of codimension~1 singularities
 with respect to the rotation of links in $V$.
\smallskip

\noindent
The \emph{third} arrow is a reformulation of the previous reduction
 in terms of canonical loops of links in the space of all links in the solid torus $V$.
\smallskip

\noindent
The \emph{fourth} arrow is a reduction of generic links to
 their 2-dimensional  diagram surfaces considered up to
 3-dimensional moves in Figure~\ref{fig:MovesDiagramSurfaces}.
\smallskip

\noindent
The \emph{fifth arrow} is a final reduction of generic links to
 their trace graphs considered up to equivalence generated by
 the moves in Figure~\ref{fig:MovesTraceGraphs}.
\smallskip

\noindent
The key ingredient of the proofs is a description of versal deformations and bifurcation diagrams
 of codimension~2 multi local singularities of plane curves in section~\ref{sect:ThroughCodim2Singularities}.
\medskip

\noindent
{\bf Acknowledgements.}
The second author is especially grateful to Hugh Morton
 and Farid Tari for fruitful comments and suggestions.
He also thanks Yu.~Burman, M.~Kazaryan, V.~Vassiliev, 
 V.~Zakalyukin for useful discussions.
\smallskip


\section{Singular subspaces in the space of all links}
\label{sect:SingularSubspaces}


\subsection{Codimension~1 singularities of link diagrams}
\noindent
\smallskip

Let $M,N$ be smooth finite dimensional manifolds.
Denote by $J^k_{[l]}(M,N)$ the space of all $l$-tuple $k$-jets of
 smooth maps $\xi:M\to N$ for all tuples $(u_1,\dots,u_l)\in M^l$,
 see \cite[sections I.2]{AVG}.
Let $(x_1,\dots,x_m)$ and $(y_1,\dots,y_n)$ be
 local coordinates in $M$ and $N$, respectively.
If the map $\xi$ is defined locally by
 $y_j=\xi_j(x_1,\dots,x_n)$, $j=1,\dots,m$, then
 the $l$-tuple $k$-\emph{jet} of the map $\xi$ at a point $(u_1,\dots,u_l)$
 is determined by
$$\{x_1,\dots,x_m\};\quad
  \{y_1,\dots,y_n\};\quad
  \left\{ \dfrac{ \bd\xi_j }{ \bd x_i } \right\};\quad
  \; \ldots \;
  \left\{ \dfrac{ \bd^{k}\xi_j }{ \bd x_{i_1}\dots \bd x_{i_s} } \right\},
  i_1+\dots +i_s=k.$$

The above quantities define local coordinates
 in $J^k_{[l]}(M,N)$.
The $l$-tuple $k$-\emph{jet} $j^k_{[l]}\xi$ of
 a smooth map $\xi:M\to N$ can be considered as
 the map $j^k_{[l]}\xi:M^l\to J^k_{[l]}(M,N)$,
 namely $(u_1,\dots,u_l)$ goes to
 the $l$-tuple $k$-jet of the map $\xi$ at $(u_1,\dots,u_l)$.
\medskip

Take an open set $W\subset J^k_{[l]}(M,N)$
 for some $k,l$.
The set of smooth maps $f:M\to N$ with $l$-tuple $k$-jets from $W$
 is called \emph{open}.
These sets for all open $W\subset J^k_{[l]}(M,N)$
 form a basis of the \emph{Whitney} topology
 in $C^{\infty}(M,N)$.
So two maps are close in the Whitney topology if
 they are close with all theirs derivatives.
\medskip

\begin{definition}
\label{def:SpaceLinks}
The \emph{space} $\SL$ of all links $K\subset V$ inherits
 the \emph{Whitney} topology from $C^{\infty}(\sqcup_{j=1}^m S_j^1,V)$.
A link $K$ defined by a smooth embedding $f:\sqcup_{j=1}^m S_j^1\to V$
 is called \emph{general} if the diagram
 $D=\pr_{xz}(K)\subset\axz$ is \emph{general}, namely
\smallskip

$\bu$
 the map $\pr_{xz}\circ f:\sqcup_{j=1}^m S_j^1\to\axz$
 is a smooth embedding outside finitely\\
 \hspace*{7mm}
 many double \emph{crossings},
 an overcrossing arc is specified at each crossing;
\smallskip

$\bu$
 the extrema of $\pr_z:D\to S_z^1$ are not crossings
 and have distinct $z$-coordinates.
\smallskip

\noindent
Denote by $\Si^{(0)}\subset\sk$ the subspace
 of all general links.
\ed
\end{definition}
\smallskip

We consider local singularities, so fix coordinates $x,z$ around each point in $\axz$.
The $x$-axis is said to be \emph{horizontal}, i.e. it is perpendicular
 to the vertical core $S_z^1\subset\axz$.
Classical codimension~1 singularities of plane curves were described by
 David \cite[List I on p.~561]{Dav}, namely the ordinary cusp $\cusp$
 (the $A_2$ singularity in Arnold's notations), simple tangency $\tang$ ($A_3$) 
 and triple point $\trip$ ($D_4$).
The solid torus $V$ has the distinguished vertical direction along $S_z^1$,
 so we also consider singularities with respect to $\pr_z:V\to S_z^1$,
 eg Reidemeister move IV is generated by passing through
 a critical crossing $\crit$, where one of the tangents is horizontal. 
\smallskip

\begin{definition}
\label{def:Codim1Singularities}
A  \emph{diagram} $D$ is the image of a smooth map $g:\sqcup_{j=1}^m S_j^1\to\axz$.
\medskip

\noindent
$\trip$
A \emph{triple point} of the diagram $D$
 is a transversal intersection $p$ of 3 arcs\\
 \hspace*{5mm}
 such that all the tangents at $p$ are not horizontal.
\medskip

\noindent
$\tang$
A \emph{simple tangency}
 is an intersection $p$ of 2 arcs given
 locally by $u=\pm v^2$.\\
 \hspace*{5mm}
We assume that the tangent at $p$ is not horizontal.
\medskip

\noindent
$\cusp$
An \emph{ordinary cusp} is the singular point $p$
 of an arc given locally by $u^2=v^3$.\\
 \hspace*{5mm}
We assume that the tangent at $p$ is not horizontal.
\medskip

\noindent
$\crit$
A \emph{critical crossing}
 is a transversal intersection $p$
 of 2 arcs such that\\
 \hspace*{5mm}
 one of the tangents at $p$ is horizontal.
\medskip

\noindent
$\cube$
A \emph{cubical point} is the singular point $p$
 of an arc given locally by $z=u^3$,\\
 \hspace*{5mm}
 the tangent at $p$ is horizontal.
\medskip

\noindent
$\maximin$
A \emph{mixed pair}
 is a pair of a local maximum and a local minimum of\\
 \hspace*{5mm}
 the projection $\pr_z:D\to S_z^1$, lying in the same horizontal line.
\medskip

\noindent
$\maximax$
An \emph{extreme pair} is a pair of
 either 2 local maxima or 2 local minima\\
 \hspace*{5mm}
 of the projection $\pr_z:D\to S_z^1$,
 lying in the same horizontal line.
\medskip

\noindent
Given a singularity
 $\ga\in\{\trip,\tang,\cusp,\crit,\cube,\maximin,\maximax\}$,
 denote by $\Si_{\ga}\subset\sk$ the singular subspace
 consisting of all links $K\subset V$ such that 
 $\pr_{xz}(K)$ is general outside $\ga$.
$$\mbox{Set }\qquad\qquad\Si^{(1)}=
  \Si_{\trip}\cup\Si_{\tang}\cup\Si_{\cusp}\cup\Si_{\crit}\cup
  \Si_{\cube}\cup\Si_{\maximin}\cup\Si_{\maximax}\subset\sk.
  \eqno{\blacksquare}$$
\end{definition}
\smallskip


\subsection{Extended Reidemeister theorem}
\label{subs:ReidemeisterTheorem}
\noindent
\smallskip

\noindent
\begin{definition}
\label{def:Transversality}
Let $M$ be a finite dimensional smooth manifold.
A subspace $\La\subset M$ is called
 \emph{a stratified space} if $\La$ is the union
 of disjoint smooth submanifolds $\La^i$ (\emph{strata})
 such that the boundary of each stratum
 is a finite union of strata of less dimensions.
Let $N$ be a finite dimensional manifold.
A smooth map $\xi:M\to N$ is \emph{transversal} to
 a smooth submanifold $U\subset N$ if
 the spaces $f_*(T_xM)$ and $T_{f(x)}U$ generate $T_{f(x)}N$
 for each $x\in M$.
A smooth map is $\eta:M\to N$ \emph{transversal} to
 a stratified space $\La\subset N$
 if the the map $\eta$ is transversal to each stratum of $\La$.
\ed
\end{definition}

\noindent
Briefly Theorem~\ref{thm:Transversality} says that any map
 can be approximated by `a nice map'.
\medskip

\noindent
\begin{theorem}
\label{thm:Transversality}
\emph{(Multi-jet transversality theorem of Thom,
 see \cite[section~I.2]{AVG})}\\
Let $M,N$ be compact smooth manifolds,
 $\La\subset J^k_{[l]}(M,N)$ be a stratified space.
Given a smooth map $\xi:M\to N$ there is
 a smooth map $\eta:M\to N$ such that
\smallskip

$\bu$
 the map $\eta$ is arbitrarily close to $\xi$
  with respect to the Whitney topology;
\smallskip

$\bu$
 the $l$-tuple $k$-jet
 $j^k_{[l]}\eta:M^l\to J^k_{[l]}(M,N)$
 is transversal to $\La\subset J^k_{[l]}(M,N)$.
\qed
\end{theorem}

\begin{lemma}
\label{lem:ComputeCodim1}
{\bf (a)}
The subspace $\Si^{(1)}$ has codimension~1 in the space $\sk$.
\smallskip

\noindent
{\bf (b)}
The subspace $\Si^{(0)}$ is open and dense in the space $\sk$.
\end{lemma}

\begin{sketch}
{\bf (a)} 
It is a standard computation in the space $J_{[3]}^1(\R,\R^2)$ 
 of 3-tuple 1-jets of maps $(x(r),z(r)):\R\to\R^2$.
For instance, fixing 3 parameters $r_1,r_2,r_3$, the subspace $\Si_{\trip}$
 maps to the subspace of $J_{[3]}^1(\R,\R^2)$ given by
 4 equations $x(r_1)=x(r_2)=x(r_3)$, $z(r_1)=z(r_2)=z(r_3)$
 and 3 inequalities $\dot z(r_i)\neq 0$, $i=1,2,3$, meaning
 that the tangents are not horizontal, hence the codimension 
 of the subspace $\Si_{\trip}\subset\SL$ is 1
 after forgetting the 3 parameters.
Analogously $\Si_{\cusp}$  maps to the subspace given by 
 4 equations $\dot x(r_1)=\dot z(r_1)=0$, $r_1=r_2=r_3$, 
 hence the codimension of $\Si_{\cusp}\subset\SL$ is 1.
A similar detailed argument will be given in the proof
 of Lemma~\ref{lem:ComputeCodim2}.
\smallskip

\noindent
{\bf (b)}
The conditions of Definition~\ref{def:SpaceLinks} define an open subset of $\SL$
 whose complement is clearly the closure of the codimension~1 subspace $\Si^{(1)}$
 from Definition~\ref{def:Codim1Singularities}.
\qed
\end{sketch}
\medskip

The following result immediately follows from Lemma~\ref{lem:ComputeCodim1} 
 since by Theorem~\ref{thm:Transversality} any isotopy in the space $\SL$ 
 of links can be approximated  by a path transversally intersecting 
 the singular subspace $\Si^{(1)}\subset\SL$ of codimension~1.
\smallskip

\begin{figure}[!h]
\includegraphics[scale=1.0]{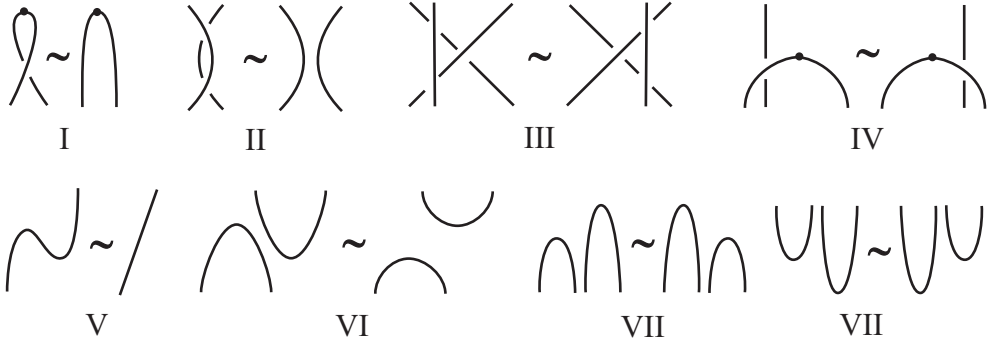}
\caption{Reidemeister moves taking into account
 local extrema}
\label{fig:ReidemeisterMoves}
\end{figure}

\begin{proposition}
\label{prop:ExtReidemeister}
Any smooth link can be approximated by a general link.
General links are \emph{isotopic} if and only if
 their diagrams can be obtained from each other by
 a plane isotopy and finitely many Reidemeister moves in 
 Figure~\ref{fig:ReidemeisterMoves}.
\qed
\end{proposition}
\smallskip

In Figure~\ref{fig:ReidemeisterMoves} orientations of arcs 
 and symmetric images of the moves are omitted.


\subsection{The co-orientation of codimension~1 subspaces}
\label{subs:Coorientation}
\noindent
\smallskip

Using Gauss diagrams of link diagrams, we define 
 the co-orientation of codimension~1 subspaces
 $\Si_{\trip},\Si_{\tang},\Si_{\cusp},\Si_{\crit}$
 from Definition~\ref{def:Codim1Singularities}.
\smallskip

\begin{definition}
\label{def:GaussDiagram}
Let a general link $K$ be defined by
 $f:\sqcup_{j=1}^m S_j^1\to V$.
The \emph{Gauss} diagram $\gd(K)$
 is the union $\sqcup_{j=1}^m S_j^1$ with chords 
 connecting points $s_1,s_2$ such that $\pr_{xz}(f(s_1))=\pr_{xz}(f(s_2))$.
Gauss diagrams $\gd_1,\gd_2$ are \emph{equivalent}
 if there is an orientation preserving diffeomorphism
 of $\sqcup_{j=1}^m S_j^1$ such that the endpoints
 of any chord of $\gd_1$ map onto the endpoints of
 a chord of $\gd_2$ and vice versa.
\ed
\end{definition}
\smallskip

\begin{figure}[]
\includegraphics[scale=1.0]{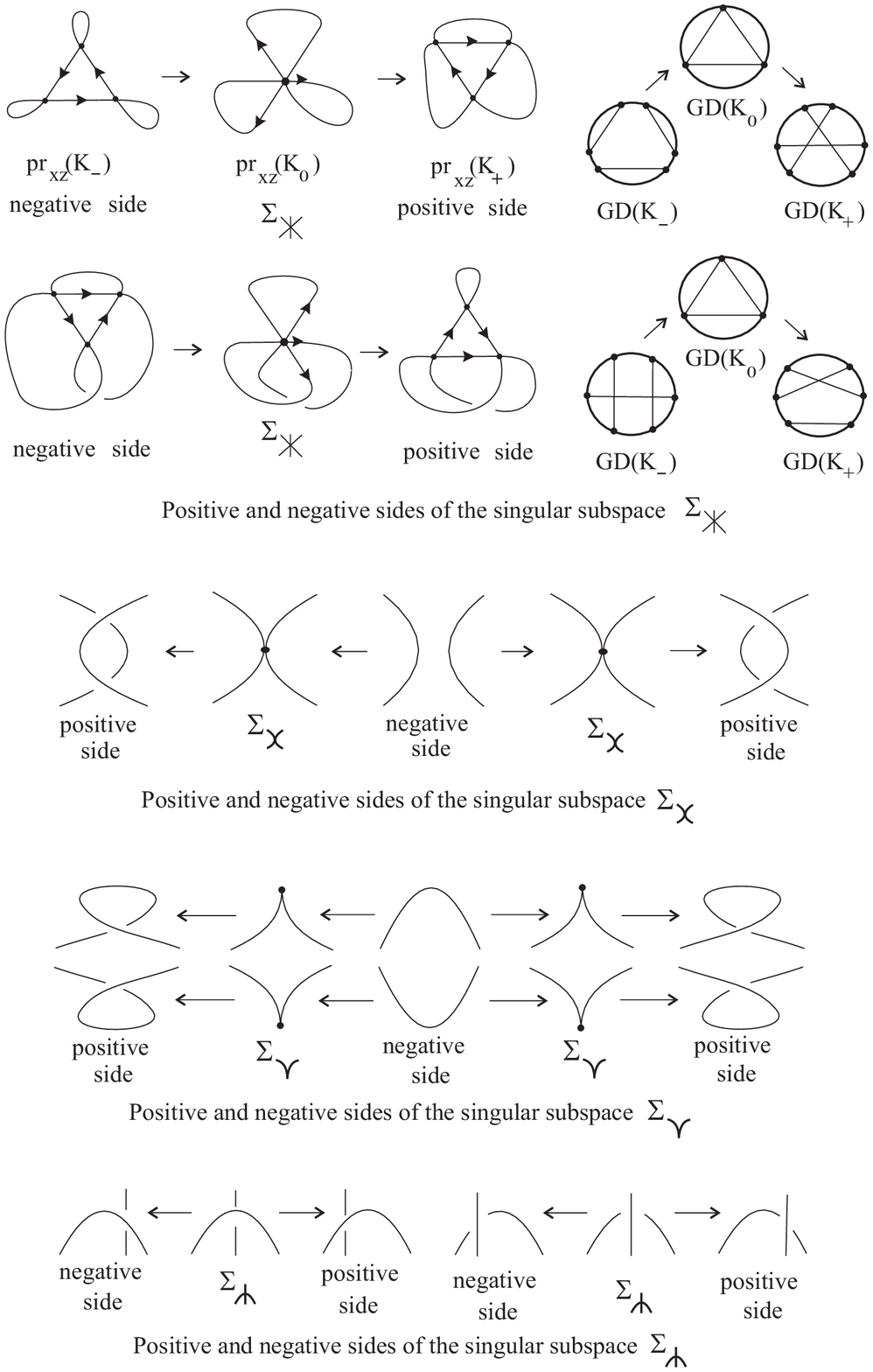}
\caption{How to define the co-orientations of 
 codimension~1 subspaces}
\label{fig:Coorientation}
\end{figure}

\begin{definition}
\label{def:Coorientation}
For each of 2 types of oriented triple points,
 the \emph{co-orientation} of $\Si_{\trip}$ is defined
 in terms of the Gauss diagrams of the corresponding links $K_{\pm}$ in 
 Figure~\ref{fig:Coorientation}.
Assume that, while $t\in S_t^1$ increases,
 the point $\rot_t(K)\in\sk$ passes through $\Si_{\trip}$
 from the negative side to the positive one.
Then associate to the corresponding triple vertex of $\tg(K)$
 the \emph{positive} sign $+$,
 otherwise take the \emph{negative} sign $-$.
The \emph{co-orientations} of $\Si_{\tang}$, $\Si_{\cusp}$,
 $\Si_{\crit}$ are similarly defined in Figure~\ref{fig:Coorientation}.
\ed
\end{definition}

Look at the trefoil $K$ in Figure~\ref{fig:RotatedTrefoils} and 
 its trace graph $\tg(K)$ in Figure~\ref{fig:TraceGraphTrefoil}.
Consider the first triple vertex of $\tg(K)$ 
 at the critical moment $t_1\in(\pi/4,\pi/2)$.
The knot $\rot_{\pi/4}(K)$
 is on the positive side of $\Si_{\trip}$ 
 (the 1st type in Figure~\ref{fig:Coorientation}), 
 while $\rot_{\pi/2}(K)$ is on the negative side of $\Si_{\trip}$, ie 
 the first triple vertex has the positive sign.
At the second triple vertex for $t_2\in(\pi/2,3\pi/4)$,
 the knot $\rot_t(K)$ goes from the negative side to the positive side.
So the second triple point also the positive sign.
\smallskip


\section{Generic links, equivalences, loops and homotopies}
\label{sect:GenericLinks}

\subsection{The canonical loop of a link and generic links}
\label{subs:CanonicalLoop}
\noindent
\smallskip

Generic links will be defined as the most generic points
 in the space $\sk$ of all links $K\subset V$ with respect to 
 the rotation $\rot_t$ of the solid torus $V$.
\smallskip

\begin{definition}
\label{def:CanonicalLoop}
The \emph{canonical} loop $\cl(K)\subset\sk$ of
 a smooth link $K\subset V$
 is the union of the rotated links $\rot_t(K)\in\sk$
 over all $t\in S_t^1$.
\smallskip

\noindent
A link $K\subset V$ is \emph{generic}
 if there are finitely many $t_1,\dots,t_k\in S_t^1$ such that
\smallskip

$\bu$ 
 for all $t\notin\{t_1,\dots,t_k\}$,
 the links $\rot_t(K)$ are general, ie $\rot_t(K)\in\Si^{(0)}$;
\smallskip

$\bu$ 
 $\cl(K)$ transversally intersects
 $\Si_{\trip}\cup\Si_{\tang}\cup\Si_{\cusp}\cup\Si_{\crit}$
 at each $t\in\{t_1,\dots,t_k\}$.
\medskip

\noindent
Denote by $\Om^{(0)}\subset\sk$ the subspace of all generic links in $V$.
\ed
\end{definition}
\smallskip

Morse modifications of index 1 would change the trace graph dramatically.
Luckily following Lemma~\ref{lem:LoopDontTouch} shows that 
 they can not occur under strong equivalence.
More exactly Lemma~\ref{lem:LoopDontTouch} shows that the canonical loop $\cl(K)$
 never touches the subspace
 $\Si_{\tang}\cup\Si_{\cusp}\cup\Si_{\crit}$
 for any link $K$.
Therefore the transversality from
 the last condition of Definition~\ref{def:CanonicalLoop}
 is relevant only for the subspace $\Si_{\trip}$.
\smallskip

\begin{lemma}
\label{lem:LoopDontTouch}
\emph{(Main topological lemma)}
For any link $K\subset V$, the canonical loop $\cl(K)$ does not touch the subspace
 $\Si_{\tang}\cup\Si_{\cusp}\cup\Si_{\crit}$.
More formally, if $K\in\Si_{\ga}$ for $\ga=\tang,\cusp,\crit$,
 the links $\rot_{\pm\e}(K)$ are on
 different sides of $\Si_{\ga}$ for small $\e>0$.
\end{lemma}
\begin{proof}
For the subspaces $\Si_{\tang}$ and $\Si_{\cusp}$,
 the projections of two small arcs with
 a tangent point (respectively, a cusp)
 are interchanged under the rotation.
\smallskip

Figure~\ref{fig:Coorientation} shows that the links $\rot_{\pm\e}(K)$
 are on different sides of $\Si_{\tang}$ and $\Si_{\cusp}$,
 respectively, since the tangent at $p$ is not horizontal, 
 ie not orthogonal to the vertical axis $S_z^1$.
The argument for $\Si_{\crit}$ is the same,
 since one tangent at the critical crossing
 is not horizontal, see the last picture of Figure~\ref{fig:Coorientation}.
\end{proof}

\begin{example}
The canonical loop $\cl(K)$ of a knot $K\subset V$
 can touch the subspace $\Si_{\trip}$.
Consider the three arcs $J_1,J_2,J_3\subset\R^3$
 defined by (see Figure~\ref{fig:TouchTriplePoints} below)
$$\left\{\begin{array}{l}
 x_1=\tau u^2,\\
 y_1=0,\\
 z_1=u;
\end{array} \right.\qquad
\left\{\begin{array}{l}
 x_2=u,\\
 y_2=-1,\\
 z_2=u;
\end{array} \right.\qquad
\left\{\begin{array}{l}
 x_3=-u,\\
 y_3=1,\\
 z_3=u;
\end{array} \right.\qquad
u\in\R,\;\tau>0.$$

Under the composition $\pr_{xz}\circ\rot_t$, 
 the arcs $J_1,J_2,J_3$ map to the following ones:
$$x_1(t)=\tau z_1^2\cos t,\quad
  x_2(t)=z_2\cos t+\sin t,\quad
  x_3(t)=-z_3\cos t-\sin t,$$
 where $z_1,z_2,z_3$ are constants.
For small $t=\e>0$, the double crossing
 $p_{23}=\pr_{xz}(\rot_{\e}(J_2))\cap\pr_{xz}(\rot_{\e}(J_3))$
 has the coordinates $x=0$, $z=-\tan\e$.
Then $p_{23}$ is at the left of the first rotated arc
 $x_1(t)=\tau z_1^2\cos t$ with respect to $\xa$.
\smallskip

For $t=-\e<0$, the crossing with $x=0$, $z=\tan\e$ is also at the left
 of the first arc.
Take a knot $K\in\Si_{\trip}$
 containing small parts of the arcs described above.
Then $\rot_{\pm\e}(K)$ are on the same side of $\Si_{\trip}$.
This means that, under the rotation of $K$, Reidemeister move~III 
 is not performed for the diagram $\pr_{xz}(\rot_t(K))$.
\end{example}

\begin{figure}[!h]
\includegraphics[scale=1.0]{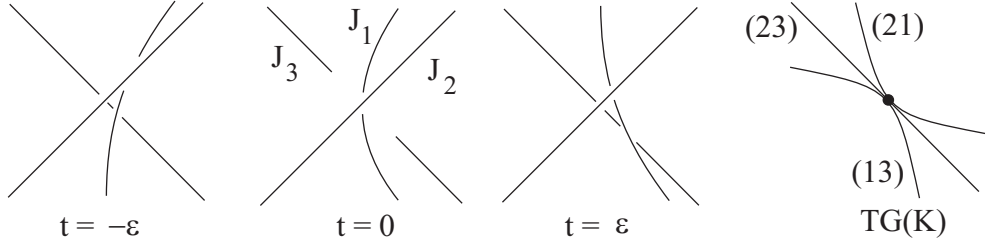}
\caption{$\cl(K)$ may touch the subspace of triple points.}
\label{fig:TouchTriplePoints}
\end{figure}


\subsection{Codimension~2 singularities and generic equivalences}
\label{subs:Codim2Singularities}
\noindent
\smallskip

Classical codimension~2 singularities of plane curves were described by
 David \cite[List II on p.~561]{Dav}, namely the rumphoidal cusp $\degcusp$
 (the $A_4$ singularity in Arnold's notations), intersected cusp $\cuspint$ ($D_5$),
 tangent triple point $\trip$ ($D_4$), cubic tangency $\cubicsn$ ($A_5$)
 and ordinary quadruple point $\quadrup$ ($X_9$).
We need to distinguish more refined singularities since the canonical loop
 of a link may not be transversal to some singular subspace, eg
 it is transversal to the codimension~2 subspace of horizontal cusps $\Si_{\horcusp}$,
 but not to the codimension~1 subspace of all cusps $\Si_{\cusp}\cup\Si_{\horcusp}$.
All tangents below are not \emph{horizontal} unless stated otherwise.
\smallskip

\begin{definition}
\label{def:Codim2Singularities}
 (\emph{codimension~2 singularities} of link diagrams)
\smallskip

\noindent
Let $D$ be a diagram, ie the image of a smooth map $g:\sqcup_{j=1}^m S_j^1\to\axz$.
\medskip

\noindent
$\quadrup$
A \emph{quadruple point} of $D$ is a transversal intersection $p$ of 4 arcs.
\medskip

\noindent
$\tangint$
A \emph{tangent triple point} of $D$
 is an intersection $p$ of 3 arcs, the first\\
 \hspace*{5mm}
 two arcs have a simple tangency and do not touch the third arc.
\medskip

\noindent
$\cuspint$
An \emph{intersected cusp} of $D$
 is an intersection of 2 arcs, where
 the first arc\\
 \hspace*{5mm}
 has an ordinary cusp whose vector $(\ddot x,\ddot z)$ 
 does not touch the second arc.
\medskip

\noindent
$\cubicsn$
A \emph{cubic tangency}
 is an intersection of 2 arcs
 given locally by $u=0$, $u=v^3$.
\medskip

\noindent
$\degcusp$
A \emph{ramphoidal cusp}
 is the singular point of an arc,
 given locally by $u^2=v^5$.
\medskip

\noindent
$\horcusp$
A \emph{horizontal cusp}
 is an ordinary cusp with horizontal tangent.
\medskip

\noindent
$\hortang$
A \emph{mixed tangency}
 is a simple tangency with a horizontal tangent
 such that\\
 \hspace*{5mm}
 one of the extrema is a maximum, another one is a minimum.
\medskip

\noindent
$\maxtang$
An \emph{extreme tangency}
 is a simple tangency with a horizontal tangent such that \\
 \hspace*{5mm}
 {\bf both} extrema are either maxima or minima.
\medskip

\noindent
$\hortrip,\inctrip$
A \emph{horizontal triple point}
 is a triple intersection, where
 the tangent line of\\
 \hspace*{5mm}
 the first arc is horizontal,
 the tangent lines of the other arcs are not horizontal.
\medskip

\noindent
Given a singularity
 $\de\in\{\quadrup,\tangint,\cuspint,\cubicsn,\degcusp,
  \hortrip,\inctrip,\horcusp,\hortang,\maxtang\}$,
 denote by $\Si_{\de}$ the union of all links $K\subset V$
 such that the diagram $\pr_{xz}(K)$ is general outside $\de$.
Set
$$\Si^{(2)}=
  \Si_{\quadrup}\cup\Si_{\tangint}\cup\Si_{\cuspint}\cup
  \Si_{\cubicsn}\cup\Si_{\degcusp}\cup\Si_{\hortrip}\cup
  \Si_{\inctrip}\cup\Si_{\horcusp}\cup\Si_{\hortang}\cup
  \Si_{\maxtang}\subset\sk.\eqno{\blacksquare}$$
\end{definition}
\smallskip

\begin{lemma}
\label{lem:ComputeCodim2}
The singular subspace $\Si^{(2)}$ has codimension 2
 in the space $\sk$.
\end{lemma}

\begin{proof}
We use multi jets of maps $(x(r),z(r)):\R\to\R^2$
 defining a diagram $D$.
Fixing 4 points $r_1,r_2,r_3,r_4$, each singularity $\de$ from 
 Definition~\ref{def:Codim2Singularities} can be described in terms of 
 4-tuple 3-jets from the space $J_{[4]}^{3}(\R,\R^2)$, where each point 
 has the 36 coordinates:
 $$J^3_{[4]}(\R,\R^2):\quad \left\{ \; r_i;\quad
  \begin{array}{cc}
   x(r_i),\quad z(r_i);\quad &
   \dot x(r_i),\quad \dot z(r_i); \\
   \ddot x(r_i),\quad \ddot z(r_i);\quad &
   \tdotx(r_i),\quad \tdotz(r_i);
  \end{array}\quad
   i=1,2,3,4. \right.$$

The jets over all $K\in\Om_{\de}$ form the finite dimensional 
 \emph{subspace} $\ti\Si_{\de}\subset J^3_{[4]}(\R,\R^2)$.
A simple tangency of 2 arcs at $r_i,r_j$ is described by 
 $\Gamma_{ij}=\deter{\dot x(r_i)}{\dot z(r_i)}{\dot x(r_j)}{\dot z(r_j)}=0.$
The frequent inequality $\dot z(r_i)\neq 0$ below says
 that the tangent of $D$ at $r_i$ is not horizontal.
The string $r_1\neq r_2\neq r_3\neq r_4$ will mean
 that $r_1,r_2,r_3,r_4$ are pairwise disjoint.
\medskip

\noindent
$\ti\Si_{\quadrup}$
$\left\{ \begin{array}{l}
  x(r_1)=x(r_2)=x(r_3)=x(r_4),\\
  z(r_1)=z(r_2)=z(r_3)=z(r_4),
  \end{array} \right.
 \begin{array}{l}
  r_1\neq r_2\neq r_3\neq r_4,\\
  \dot z(r_i)\neq 0,\; \Gamma_{ij}\neq 0,\; i\neq j;
 \end{array}$
\medskip

\noindent
$\ti\Si_{\tangint}$
$\left\{ \begin{array}{l}
  x(r_1)=x(r_2)=x(r_3),\\
  z(r_1)=z(r_2)=z(r_3),\\
   \end{array} \right.
 \begin{array}{l}
  r_1\neq r_2\neq r_3=r_4,\; \dot z(r_i)\neq 0,\\
  \Gamma_{12}=0,\; \Ga_{23}\neq 0,\; \Ga_{13}\neq 0;
 \end{array}$
\medskip

\noindent
$\ti\Si_{\cuspint}$
$\left\{ \begin{array}{l}
 x(r_1)=x(r_2),\; z(r_1)=z(r_2), \\
 \dot x(r_1)=\dot z(r_1)=0, 
 \end{array} \right.$
$\begin{array}{l}
 \ddot z(r_1)\neq 0,\; \dot z(r_2)\neq 0, \\
 r_1\neq r_2=r_3=r_4,
 \end{array} 
\deter{\ddot x(r_1)}{\ddot z(r_1)}{\dot x(r_2)}{\dot z(r_2)}\neq 0;$
\medskip

\noindent
$\ti\Si_{\cubicsn}$
$\left\{ \begin{array}{l}
 x(r_1)=x(r_2),\; z(r_1)=z(r_2),\\
   r_1\neq r_2=r_3=r_4,\; \Ga_{12}=0, \\
   \dot z(r_1) \neq 0, \quad \dot z(r_2) \neq 0, \\
  \end{array} \right. \;
 \deter{\ddot x(r_1)}{\ddot z(r_1)}{\dot x(r_2)}{\dot z(r_2)}=0,\;
 \deter{\tdotx(r_1)}{\tdotz(r_1)}{\dot x(r_2)}{\dot z(r_2)}\neq 0;$
\medskip

\noindent
$\ti\Si_{\degcusp}$
$\left\{
 \begin{array}{c}
  \dot x(r_1)=\dot z(r_1)=0,\quad \ddot z(r_1)\neq 0\\
  r_1=r_2=r_3=r_4,
 \end{array} \right. \quad
 \left| \begin{array}{cc}
  \ddot x(r_1) & \tdotx(r_1) \\
  \ddot z(r_1) & \tdotz(r_1) \\
 \end{array} \right| =0;$
\medskip

\noindent
$\ti\Si_{\hortrip}\cup \ti\Si_{\inctrip}$
$\left\{ \begin{array}{l}
  x(r_1)=x(r_2)=x(r_3),\\  
    z(r_1)=z(r_2)=z(r_3),
 \end{array} \right.
 \begin{array}{l}
  \dot z(r_1)=0,\;
  \ddot z(r_1)\neq 0, \\
  \dot z(r_2)\neq 0,\;
  \dot z(r_3)\neq 0,
 \end{array}
 \begin{array}{l}
    r_1\neq r_2\neq r_3=r_4,\\
    \Ga_{ij}\neq 0,\; i\neq j;
 \end{array}$
\medskip

\noindent
$\quad\ti\Si_{\horcusp}$
$\quad \{ \quad \dot x(r_1)=\dot z(r_1)=\ddot z(r_1)=0$,
$\quad r_1=r_2=r_3=r_4$,
$\quad \tdotz(r_1)\neq 0$;
\medskip

\noindent
$\ti\Si_{\hortang}\cup \ti\Si_{\maxtang}$
$\left\{ \begin{array}{l}
  z(r_1)=z(r_2),\quad \dot z(r_1)=\dot z(r_2)=0, \\
  x(r_1)=x(r_2), \quad r_1\neq r_2=r_3=r_4,
 \end{array} \right.
 \begin{array}{l}
 \ddot z(r_1)\neq 0,\\
 \ddot z(r_2)\neq 0;
 \end{array}$
\medskip

The conditions above can be obtained using classical normal forms of the singularities,
 eg the ramphoidal cusp $\degcusp$ is a degeneration of the ordinary cusp $\cusp$ 
 clearly given by $\dot x(r)=\dot z(r)=0$.
Locally one has $(x,z)=(a_2r^2+a_3r^3+\dots,b_2r^2+b_3r^3+\dots)$, $b_2\neq 0$, 
 which is (left) equivalent to $(x,z)=((a_3-b_3 a_2/b_2)r^3+\dots,b_2 r^2+\dots)$, hence
 $a_2b_3=b_2a_3$, ie the vectors $(\ddot x(r),\ddot z(r))$ and $(\tdotx(r),\tdotz(r))$ are collinear.
\smallskip

Each subspace $\ti\Si_{\de}$ is defined by 6 equations,
 hence $\codim\ti\Om_{\de}=6$ in $J^3_{[4]}(\R,\R^3)$.
The subspaces $\Si_{\de}$ from Definition~\ref{def:Codim2Singularities} 
 map to the corresponding subspaces $\ti\Si_{\de}\subset J^3_{[4]}(\R,\R^2)$
 by adding 4 points $r_1,r_2,r_3,r_4$ on a diagram.
When we forget these points, the codimension decreases by 4,
 ie $\codim\Si_{\de}=2$ in the space $\sk$.
\end{proof}
\smallskip

\begin{definition}
\label{def:GenericEquivalence}
Let $\Om_{\tantrip}$ be the set of all links failing to be generic
 due to exactly one tangency of $\cl(K)$ with 
 the codimension~1 subspace $\Si_{\trip}$.
\smallskip

\noindent
Given $\de\in\{\quadrup,\tangint,\cuspint,\cubicsn,\degcusp,
  \hortrip,\inctrip,\horcusp,\hortang,\maxtang\}$,
 let $\Om_{\de}$ consist of all links $K$ failing 
 to be generic because of exactly one transversal intersection
 of $\cl(K)$ with $\Si_{\de}$.
Set
$$\Om^{(1)}=
  \Om_{\quadrup}\cup\Om_{\tangint}\cup\Om_{\cuspint}\cup
  \Om_{\cubicsn}\cup\Om_{\degcusp}\cup\Om_{\hortrip}\cup
  \Om_{\inctrip}\cup\Om_{\horcusp}\cup\Om_{\hortang}\cup
  \Om_{\maxtang}\cup\Om_{\tantrip}.$$

\noindent
A \emph{generic equivalence} is a smooth path
 $F:[0,1]\to\sk$ intersecting transversally
 the subspace $\Om^{(1)}$, ie
 there are finitely many $r_1,\dots,r_k\in[0,1]$ such that
\smallskip

$\bu$
 the links $F(r)\in\sk$ are generic for all
 $r\notin\{r_1,\dots,r_k\}$;
\smallskip

$\bu$
 the canonical loop $\cl(F(r))$ transversally intersects
 $\Om^{(1)}$ for $r=r_1,\dots,r_k$.
\ed
\end{definition}
\smallskip

\begin{lemma}
\label{lem:ComputeCodim}
{\bf (a)}
The subspace $\Om^{(1)}$ has codimension~1 in the space $\sk$.
\smallskip

\noindent
{\bf (b)}
The subspace $\Om^{(0)}$ is open and dense in the space $\sk$.
\end{lemma}

\begin{proof}
{\bf (a)}
Choose a link $K\subset V$ 
 given by an embedding $f:\sqcup_{j=1}^m S_j^1\to V$
 that fails to be generic due to exactly one singularity $\de$ from Definition~\ref{def:GenericEquivalence}.
These singularities were introduced using the rotation of the solid torus $V$.
So we describe them in terms of maps $\R\to\R^3$,
 not $\R\to\R^2$ as in the proof of Lemma~\ref{lem:ComputeCodim2}.
\smallskip

There is a 4-tuple $(r_1,r_2,r_3,r_4)\in (\sqcup_{j=1}^m S_j^1)^4$
 defining the chosen singularity of $K\in\Om_{\de}$.
The 4-tuple 3-jet $j^3_{[4]}f(r_1,r_2,r_3,r_4)$
 is a point in $J^3_{[4]}(\R,\R^3)$.
These points over all $K\in\Om_{\de}$ form the finite dimensional 
 \emph{subspace} $\ti\Om_{\de}\subset J^3_{[4]}(\R,\R^3)$.

We check that $\ti\Om_{\de}$ has codimension~5 in
 $J^3_{[4]}(\R,\R^3)$.
Denote by $x(r),y(r),z(r)$ the compositions of 
 $f:\sqcup_{j=1}^m S_j^1\to V\subset\R^3$ and the projections 
 to the coordinate axes.
Then the 4-tuple 3-jet of $K$ is determined by the following 52 quantities.
 $$J^3_{[4]}(\R,\R^3):\; \left\{ \; r_i;\quad
  \begin{array}{cc}
   x(r_i),\quad y(r_i),\quad z(r_i);\quad &
   \dot x(r_i),\quad \dot y(r_i),\quad \dot z(r_i); \\
   \ddot x(r_i),\quad \ddot y(r_i),\quad \ddot z(r_i);\quad &
   \tdotx(r_i),\quad \tdoty(r_i),\quad \tdotz(r_i);
  \end{array}\;
   i=1,2,3,4. \right.$$
For $i,j\in\{1,2,3,4\}$, $i\neq j$,
 introduce the differences
$$\De x_{ij}=x(r_i)-x(r_j),\quad
  \De y_{ij}=y(r_i)-y(r_j),\quad
  \De z_{ij}=z(r_i)-z(r_j).$$

Points $f(r_i),f(r_j),f(r_k)\in K$
 project to the same point under $\pr_{xz}:\rot_t(K)\to\axz\times\{t\}$
 for some $t$ if and only if $z(r_i)=z(r_j)=z(r_k)$,
 $\deter{\De x_{ij}}{\De y_{ij}}{\De x_{jk}}{\De y_{jk}}=0.$
The last determinant is (up to the sign) the area of the triangle with
 the vertices $(x(r_i),y(r_i))$, $(x(r_j),y(r_j))$, $(x(r_k),y(r_k))$
 in the horizontal plane $\{z(r_i)=z(r_j)=z(r_k)\}$.
\smallskip

Set $\De_{ij}=
  \left| \begin{array}{ccc}
   \dot x(r_i) & \dot x(r_j) & \De x_{ij} \\
   \dot y(r_i) & \dot y(r_j) & \De y_{ij} \\
   \dot z(r_i) & \dot z(r_j) & \De z_{ij}
  \end{array} \right|.$
The diagram $\pr_{xz}(\rot_t(K))$ contains
 two arcs having a simple tangency at $r=r_i$, $r=r_j$ and some $t$
 if and only if $z(r_i)=z(r_j)$ and $\De_{ij}=0$, ie
 the straight line through $f(r_i),f(r_j)\in K$ lies in the plane
 spanned by the tangent vectors of $K$ at $r=r_i$ and $r=r_j$.
\smallskip

We describe analytically the subspaces $\ti\Om_{\de}$
 associated to the singularities
$$\quadrup,\tangint,\cuspint,\cubicsn,\degcusp,
 \hortrip,\inctrip,\hortang,\horcusp,\maxtang,\tantrip.$$

\noindent
$\ti\Om_{\quadrup}$
$\left\{ \begin{array}{l}
  z(r_1)=z(r_2)=z(r_3)=z(r_4),\\
  r_1\neq r_2\neq r_3\neq r_4,\\
  \dot z(r_i)\neq 0,\; i=1,2,3,4,
  \end{array} \right.
 \begin{array}{l}
  \deter{\De x_{12}}{\De y_{12}}{\De x_{23}}{\De y_{23}}=
  \deter{\De x_{12}}{\De y_{12}}{\De x_{24}}{\De y_{24}}=0,\\
  \De_{ij}\neq 0,\; i,j\in\{1,2,3,4\},\; i\neq j;
 \end{array}$
\medskip

\noindent
$\ti\Om_{\tangint}$
$\left\{ \begin{array}{l}
  z(r_1)=z(r_2)=z(r_3),\\
  r_1\neq r_2\neq r_3=r_4,\\
  \dot z(r_i)\neq 0,\; i=1,2,3,
 \end{array} \right.
 \begin{array}{l}
  \deter{\De x_{12}}{\De y_{12}}{\De x_{23}}{\De y_{23}}=0,\\
  \De_{12}=0,\; \De_{23}\neq 0,\; \De_{13}\neq 0;
 \end{array}$
\medskip

\noindent
$\ti\Om_{\cuspint}$
$\left\{ \begin{array}{l}
 z(r_1)=z(r_2), \\
 \dot z(r_1)=0,\; \ddot z(r_1)\neq 0, \\
 \dot z(r_2)\neq 0,
\end{array} \right.$
$\begin{array}{l}
 \deter{\dot x(r_1)}{\dot y(r_1)}{\dot x(r_2)}{\dot y(r_2)}=0,\\
 r_1\neq r_2=r_3=r_4,
 \end{array}$
 $\left| \begin{array}{ccc}
  \ddot x(r_1) & \dot x(r_2) & \De x_{12} \\
  \ddot y(r_1) & \dot y(r_2) & \De y_{12} \\
  \ddot z(r_1) & \dot z(r_2) & \De z_{12}
 \end{array} \right| \neq 0$;
\medskip

\noindent
$\ti\Om_{\cubicsn}$
$\left\{ \begin{array}{l}
   z(r_1)=z(r_2),\quad \De_{12}=0,\\
   r_1\neq r_2=r_3=r_4, \\
   \dot z(r_1) \neq 0, \quad \dot z(r_2) \neq 0, \\
  \end{array} \right.
  \left| \begin{array}{ccc}
   \ddot x(r_1) & \ddot x(r_2) & \De x_{12} \\
   \ddot y(r_1) & \ddot y(r_2) & \De y_{12} \\
   \ddot z(r_1) & \ddot z(r_2) & \De z_{12} \\
  \end{array} \right| =0;$
\medskip

\noindent
$\ti\Om_{\degcusp}$
$\left\{
 \begin{array}{c}
  \dot z(r_1)=0,\quad \ddot z(r_1)\neq 0,\\
  r_1=r_2=r_3=r_4,
 \end{array} \right. \quad
 \dfrac{\tdotx(r_1)}{\ddot x(r_1)}=
 \dfrac{\tdoty(r_1)}{\ddot y(r_1)}=
 \dfrac{\tdotz(r_1)}{\ddot z(r_1)}.$
\medskip

The last equations with 3 fractions mean that the vectors
 of the 2nd and 3rd derivatives are collinear, which
 corresponds to the similar condition for $\ti\Si_{\degcusp}$
 in the proof of Lemma~\ref{lem:ComputeCodim2}.
If a denominator is zero, the numerator must be also zero.
\medskip

\noindent
$\ti\Om_{\hortrip}\cup \ti\Om_{\inctrip}$
$\left\{ \begin{array}{l}
  z(r_1)=z(r_2)=z(r_3),\quad r_1\neq r_2\neq r_3=r_4,\\
  \dot z(r_1)=0,\;
  \ddot z(r_1)\neq 0,\;
  \dot z(r_2)\neq 0,\;
  \dot z(r_3)\neq 0,
 \end{array} \right.\;
 \deter{\De x_{12}}{\De y_{12}}{\De x_{23}}{\De y_{23}}=0$;
\medskip

\noindent
$\quad\ti\Om_{\horcusp}$
$\quad \{ \quad \dot z(r_1)=\ddot z(r_1)=0$,
$\quad r_1=r_2=r_3=r_4$,
$\quad \tdotz(r_1)\neq 0$;
\medskip

\noindent
$\ti\Om_{\hortang}\cup \ti\Om_{\maxtang}$
$\left\{ \begin{array}{l}
  z(r_1)=z(r_2),\\ 
  \dot z(r_1)=\dot z(r_2)=0, 
 \end{array} \right.
 \begin{array}{l}
 r_1\neq r_2=r_3=r_4, \\
 \ddot z(r_1)\neq 0,\; \ddot z(r_2)\neq 0.
 \end{array}$
\medskip

If $\dot z(r_i)\neq 0$, then locally $r_i$ can be considered as a function of $z$,
 hence any function of (several) $r_i$ can be differentiated with respect to $z$.
Below the tangency with $\Si_{\trip}$ means that the derivative of 
 the vanishing determinant $\De=\deter{\De x_{12}}{\De y_{12}}{\De x_{23}}{\De y_{23}}$
 defining a triple point under the projection $\pr_{xz}:\rot_t(K)\to\axz\times\{t\}$ also vanishes.
\medskip

$\ti\Om_{\tantrip}$
$\left\{ \begin{array}{l}
  z(r_1)=z(r_2)=z(r_3),\quad r_1\neq r_2\neq r_3=r_4\\
  \De=\dfrac{d}{dz}\De=0, \; \dfrac{d^2}{dz^2}\De\neq 0, \; 
  \De=\deter{\De x_{12}}{\De y_{12}}{\De x_{23}}{\De y_{23}}
 \end{array} 
 \begin{array}{l}
  \De_{ij}\neq 0,\\ 
  i\neq j,\\
  \dot z(r_i)\neq 0.
  \end{array} \right.$
\medskip

\noindent
Generic inequalities $\dfrac{dg}{dz}\neq 0$ should be added to the descriptions above
 for each condition $g=0$, which guarantees no tangency of the canonical loop
 with the corresponding subspace $\Si_{\de}$.
In important cases like $\ti\Om_{\degcusp}$ we explicitly accompanied $\dot z(r_1)=0$
 with $\ddot z(r_2)\neq 0$ equivalent to $\dfrac{\dot z(r_1(z))}{dz}\neq 0$, 
 but also every equation like $z(r_1)=z(r_2)$ should be accompanied 
 with $\dfrac{dz(r_1(z))}{dz}\neq \dfrac{dz(r_2(z))}{dz}$.
\medskip

Each subspace $\ti\Om_{\de}$ is defined by 5 equations,
 hence $\codim\ti\Om_{\de}=5$ in $J^3_{[4]}(\R,\R^3)$.
The subspaces $\Om_{\de}$ introduced geometrically
 in Definition~\ref{def:GenericEquivalence} correspond to $\ti\Om_{\de}$ 
 by adding 4 points $r_1,r_2,r_3,r_4$ on a link.
When we forget about these points the codimension decreases by 4,
 ie $\codim\Om_{\de}=1$ in the space $\sk$ of all links $K\subset V$.
\medskip

\noindent
{\bf (b)}
The conditions of Definition~\ref{def:GenericEquivalence} define an open subset of $\SL$
 whose complement is clearly the closure of the codimension~1 subspace $\Om^{(1)}$.
\end{proof}

The following result similar to Proposition~\ref{prop:ExtReidemeister}
 follows from Lemma~\ref{lem:ComputeCodim} since by Theorem~\ref{thm:Transversality} 
 any isotopy in the space $\SL$ of links can be approximated by a path 
 transversally intersecting the singular subspace $\Om^{(1)}\subset\SL$.
\smallskip

\begin{proposition}
\label{prop:ReductionToGenericLinks}
{\bf (a)}
Any smooth link can be approximated by a generic link.
\smallskip

\noindent
{\bf (b)}
Any smooth equivalence of links can be approximated
 by a generic one.
\qed
\end{proposition}
\smallskip


\subsection{Generic loops and generic homotopies in the space of links}
\label{subs:GenericLoopsHomotopies}
\noindent
\smallskip

A \emph{loop} of links $\{K_t\}\subset\sk$
 means a \emph{smooth} loop, ie
 a smooth map $S_t^1\to\sk$.
Generic loops provide a suitable generalization of
 the canonical loop.

\begin{definition}
\label{def:GenericLoop}
A smooth loop of links $\{K_t\}\subset\sk$, $t\in S_t^1$, is called \emph{generic}
 if there are finitely many critical moments
 $t_1,\dots,t_k\in S_t^1$ such that
\smallskip

$\bu$
the link $K_t$ maps to $K_{t+\pi}$ under the rotation
 through $\pi$ for every $t\in S_t^1$;
\smallskip

$\bu$
for all $t\notin\{t_1,\dots,t_k\}$,
 the links $K_t$ are general, ie $K_t\in\Si^{(0)}$;
\smallskip

$\bu$
$\{K_t\}$ transversally intersects
 $\Si_{\trip}\cup\Si_{\tang}\cup\Si_{\cusp}\cup\Si_{\crit}$
 at each $t=t_1,\dots,t_k$.
\ed
\end{definition}
\smallskip

Due to Lemmas~\ref{lem:ComputeCodim1}, \ref{lem:ComputeCodim} 
 any loop can be approximated by a generic loop.
But a generic loop may be too trivial.
For instance, a loop $S_t^1\to\sk$ contractible
 to a generic link through generic links
 carries information about only one diagram.
More interesting objects are generic loops
 homotopic to canonical loops.
\smallskip

\begin{definition}
\label{def:GenericHomotopy}
A smooth family $\{L_s\}$ of loops, $s\in[0,1]$,
 is called a \emph{generic} homotopy
 if there are finitely many
 critical moments $s_1,\dots,s_k\in[0,1]$ such that
\smallskip

\noindent
$\bu$
for $s\notin\{s_1,\dots,s_k\}$,
 the loop $L_s$ is generic
 in the sense of Definition~\ref{def:GenericLoop};
\smallskip

\noindent
$\bu$
for each $s\in\{s_1,\dots,s_k\}$,
 the loop $L_s$ fails to be generic since either
 $L_s$ transversally intersects $\Si^{(2)}$
 or $L_s$ touches $\Si_{\trip}$ at exactly one point.
\ed
\end{definition}

\begin{lemma}
\label{lem:HomotopyLoopsLinks}
{\bf (a)}
The canonical loop of any generic link is a generic loop.
\smallskip

\noindent
{\bf (b)}
Any generic equivalence $\{K_s\}$, $s\in[0,1]$, of links
 provides the generic homotopy of loops $\{\cl(K_s)\}$ of links.
\smallskip

\noindent
{\bf (c)}
If canonical loops $\cl(K_0)$ and $\cl(K_1)$
 of generic links $K_0$ and $K_1$ are generically homotopic
 then $K_0$ and $K_1$ are generically equivalent.
\end{lemma}
\begin{proof}
{\bf (a)}
The canonical loop of any link is symmetric in the sense
 that $\rot_t(K)$ maps to $\rot_{t+\pi}(K)$
 under the rotation through $\pi$ for every $t\in S_t^1$.
The other conditions of Definition~\ref{def:GenericLoop}
 correspond to the conditions of Definition~\ref{def:CanonicalLoop}.
\smallskip

\noindent
{\bf (b)}
Compare Definition~\ref{def:GenericEquivalence}
 with Definitions~\ref{def:GenericLoop} and \ref{def:GenericHomotopy}.
\smallskip

\noindent
{\bf (c)}
Let $\{L_s\}$, $s\in[0,1]$, be a generic homotopy
 between $\cl(K_0)$ and $\cl(K_1)$.
The loops $L_s$ can be represented by
 a cylinder $S_t^1\times[0,1]$
 mapped to the space $\sk$.
Take a smooth path connecting $K_0$ and $K_1$
 inside the cylinder.
This smooth equivalence can be approximated
 by a generic one due to Proposition~\ref{prop:ReductionToGenericLinks}b.
\end{proof}
\smallskip

By Lemma~\ref{lem:HomotopyLoopsLinks} the classification of links reduces
 to their canonical loops.
\smallskip

\begin{proposition}
\label{prop:ReductionToGenericLoops}
Generic links are generically equivalent in $V$ if and only if
 their canonical loops are generically homotopic
 in the space $\sk$ of all links $K\subset V$.
\qed
\end{proposition}
\smallskip


\section{Through codimension~2 singularities}
\label{sect:ThroughCodim2Singularities}


\subsection{Versal deformations of codimension~2 singularities}
\label{subs:VersalDeformations}
\noindent
\smallskip

To understand what happens when the canonical loop of a link 
 passes through the singular subspace $\Si^{(2)}$, we study 
 bifurcation diagrams of codimension~2 singularities.

\begin{lemma}
\label{lem:NormalForms}
The codimension~2 singularities from Definition~\ref{def:Codim2Singularities}
 have the normal forms in the table below, where $r$ is the parameter on the curve and
\smallskip

\noindent
$\bu$
$\Ae$ is the extended right-left equivalence, 
 i.e. diffeomorphisms of $\R^2$ don't fix 0;
\smallskip

\noindent
$\bu$
 $\Az$ is the restricted right-left equivalence such that 
 left  diffeomorphisms of $\R^2$ have the form $(g(x,z),h(z))$, 
 where $g(x,z):\R^2\to\R$, $h(z):\R\to\R$ are diffeomorphisms.
\medskip

\begin{tabular}{c|l}
$\quadrup$, $\Ae$  & $\{ x=0, \, z=r \},\; \{ x=r, \, z=r \},\; \{x=-r, \, z=r \},\; \{ x=er, \, z=r \}$\\
\hline

$\tangint$, $\Ae$ & $\{x=r^2, \, z=r\},\; \{x=0, \, z=r\},\; \{x=r, \, z=r\}$\\
\hline

$\cuspint$, $\Ae$ & $\{x=r^3,\, z=r^2\},\; \{x=r, \, z=r\}$\\
\hline

$\cubicsn$, $\Ae$ & $\{x=r^3,\, z=r\},\; \{x=0, \, z=r\}$\\
\hline

$\degcusp$, $\Ae$ & $\{x=r^5,\, z=r^2\}$\\
\hline

$\horcusp$, $\Az$ & $\{x=r^2,\, z=r^3\}$\\
\hline

$\hortang$, $\Az$ & $\{x=r,\, z=r^2\},\; \{x=r,\, z=-r^2\}$\\
\hline

$\maxtang$, $\Az$ & $\{x=r,\, z=-2r^2\},\; \{x=r,\, z=-r^2\}$\\
\hline

$\hortrip$, $\Az$ & $\{ x=r, \, z=-r^2 \},\; \{x=r, \, z=r\},\; \{x=-r, \, z=r\}$\\
\hline

$\inctrip$, $\Az$ & $\{ x=r, \, z=-r^2 \},\; \{x=r, \, z=r\},\; \{x=2r, \, z=r\}$\\
\hline
\end{tabular}
\end{lemma}
\smallskip

\begin{sketch}
The normal forms up to $\Ae$-equivalence are classical, eg
 the parameter $e\neq 0$ in the normal form of $\quadrup$ ($X_9$) can not be skipped
 as the cross-ratio of 4 slopes is invariant under diffeomorphisms,  
 see \cite[Lemma~6.5]{Wal}.
The singularities $\hortang$, $\horcusp$, $\maxtang$, $\hortrip$, $\inctrip$
 should be considered up to $\Az$-equivalence respecting $\{z=\const\}$, 
 otherwise they don't have codimension~2, e.g. the normal form $(r^2,r^3)$ of $\horcusp$
 is not $\Az$-equivalent to the normal form $(r^3,r^2)$ of $\cusp$. 
Deducing new normal forms is similar, eg the horizontal cusp $\horcusp$
 is defined by the conditions $\dot x(0)=\dot z(0)=\ddot z(0)=0$, hence  
 $x(r)=ar^2+\dots$, $z(r)=br^3+\dots$, which normalises to $(r^2,r^3)$ as required.
\qed
\end{sketch}
\medskip

Mancini and Ruas \cite{MR} have shown that the group $\Az$ from Lemma~\ref{lem:NormalForms} 
 is geometric in the sense of Damon \cite{Dam}.
So the standard technique of singularity theory can be applied to find 
 versal deformations of corresponding codimension~2 singularities.
\smallskip

We consider horizontal triple points $\hortrip$ and $\inctrip$ separately,
 because the associated moves on trace graphs look slightly different 
 in Figure~\ref{fig:MovesTraceGraphs}ix, \ref{fig:MovesTraceGraphs}x.
A deformation of a germ $(x(r),z(r)):\R\to\R^2$ with parameters $a,b$ 
 is a germ $F:\R\times\R^2\to\R^2$ such that $F(r;0,0)\equiv (x(r),z(r))$.
A deformation $F$ is \emph{versal} if any other deformation can be obtained from 
 $F$ by actions of the corresponding group $\Ae$ or $\Az$.
\smallskip

The versality can be checked using the following tangent spaces 
 at germs in the space of deformations.
Let $T^r$ be the \emph{right} tangent space at a germ $(x(r),z(r))$ 
 generated by the right diffeomorphisms $\R\to\R$, eg the right space
 $T^r$ at $(r^5,r^2)$ of $\degcusp$ consists of $(5r^4f(r),2rf(r))$,
 where $f:\R\to\R$.
Denote by $T^l$ the \emph{left} tangent space at a germ $(x(r),z(r))$ 
 generated by the restricted left diffeomorphisms $(g(x,z),h(z)):\R^2\to\R^2$, 
 where $g:\R^2\to\R$, $h:\R\to\R$ are any germs.
For instance, the left space $T^l$ at $(r^2,r^3)$ of $\horcusp$ is formed by 
 $(g(r^2,r^3),h(r^3))=(a_1+a_2r^2+a_3r^3+\dots,b_1+b_2r^3+\dots)$.
The \emph{parameter} normal space $N^p$ of a deformation $F(r;a,b)$
 consists of linear combinations $c\dfrac{\bd F}{\bd a}+d\dfrac{\bd F}{\bd b}$ 
 at $a=b=0$, where $c,d$ are constants, e.g. the space $N^p$ 
 of $(r^5+ar^3+br,r^2)$ consists of vectors $(cr^3+dr,0)$.
\smallskip

In the case of a multi-germ the right space $T_i^r$ is associated to 
 the independent right diffeomorphisms $f_i(r)$ around each point $r_i$.
The left space $T_i^l$ is generated by the same left diffeomorphisms at every $r_i$,
The parameter space $N_i^p$ is spanned by the derivatives 
 along the parameters of the deformation at each $r_i$.
\smallskip

The following standard statement is a simple application of \cite[section~I.8.2]{AVG}.

\begin{proposition}
\label{prop:CriterionVersality}
A deformation $F(r;a,b)$ of a multi-germ $(x(r),z(r)):\R\to\R^2$ is \emph{versal} if
 at every point $r_i$ any germ $\R\to\R^2$ can be represented as 
 a sum of vectors from the spaces $T_i^r$, $T_i^l$ and $N_i^p$. 
\qed
\end{proposition}

\begin{lemma}
\label{lem:VersalDeformations}
The codimension~2 singularities from Definition~\ref{def:Codim2Singularities}
 have the versal deformations with parameters $a,b$ in the table below.
\end{lemma}
\medskip

\hspace*{-5mm}
\begin{tabular}{c|l}
$\quadrup$, $\Ae$  & $\{ x=0, z=r \},\, \{ x=r+a, z=r \},\, \{x=-r-b, z=r\},\, \{ x=er, z=r  \}$\\
\hline

$\tangint$, $\Ae$ & $\{x=r ^2-2a, z=r \},\; \{x=0, z=r \},\; \{x=r-b, z=r \}$\\
\hline

$\cuspint$, $\Ae$ & $\{x=r^3-br,\, z=r^2\},\; \{x=r-a, z=r \}$\\
\hline

$\cubicsn$, $\Ae$ & $\{x=r^3-3br+a,\, z=r\},\; \{x=0, z=r \}$\\
\hline

$\degcusp$, $\Ae$ & $\{x=r^5+ar^3+br,\, z=r^2\}$\\
\hline

$\horcusp$, $\Az$ & $\{x=r^2,\, z=r^3+ar^2-br\}$\\
\hline

$\hortang$, $\Az$ & $\{x=r, z=r^2-b\},\; \{x=r+a, z=-r^2 \}$\\
\hline

$\maxtang$, $\Az$ & $\{x=r, z=-2r^2-b \},\; \{x=r+a, z=-r^2\}$\\
\hline

$\hortrip$, $\Az$ & $\{ x=r, \,z=-r^2 \},\; \{x=r+a, z=r \},\; \{x=-r-b, z=r \}$\\
\hline

$\inctrip$, $\Az$ & $\{x=r, z=-r^2\},\; \{x=r+a, z=r \},\; \{x=r/2-b, z=r \}$\\
\hline

\end{tabular}
\smallskip

\begin{sketch}
Versal deformations of classical codimension~2 singularities 
 $\degcusp$ ($A_4$), $\cuspint$ ($D_5$), $\tangint$ ($D_4$), $\cubicsn$ ($A_5$)
 and $\quadrup$ ($X_9$) up to $\Ae$-equivalence 
 were recently described by Wall \cite[subsection~6.1]{Wal}.
The remaining cases follow from the table below.
\medskip

\begin{tabular}{c|l|l|l}
singularity & $T_i^r$ & $T_i^l$ & $N_i^p$ \\
\hline

\horcusp & $(2rf(r),3r^2f(r))$ & $(g(r^2,r^3),h(r^3))$ & $(0,cr^2-dr)$\\
\hline

\hortang & $(f_1(r),2rf_1(r))$ & $(g(r,r^2),h(r^2))$ & $(0,-d)$\\
 & $(f_2(r),-2rf_2(r))$ & $(g(r,-r^2),h(-r^2))$ & $(c,0)$ \\
\hline

\maxtang & $(f_1(r),-4rf_1(r))$ & $(g(r,-2r^2),h(-2r^2))$ & $(0,-d)$\\
 & $(f_2(r),-2rf_2(r))$ & $(g(r,-r^2),h(-r^2))$ & $(c,0)$ \\
\hline

\hortrip & $(f_1(r),-2rf_1(r))$ & $(g(r,-r^2),h(-r^2))$ & $(0,0)$\\
 & $(f_2(r),f_2(r))$ & $(g(r,r),h(r))$ & $(c,0)$\\
 & $(-f_3(r),f_3(r))$ & $(g(-r,r),h(r))$ & $(-d,0)$\\
\hline

\inctrip & $(f_1(r),-2rf_1(r))$ & $(g(r,-r^2),h(-r^2))$ & $(0,0)$\\
 & $(f_2(r),f_2(r))$ & $(g(r,r),h(r))$ & $(c,0)$\\
 & $(-f_3(r)/2,f_3(r))$ & $(g(-r/2,r),h(r))$ & $(-d,0)$\\
\hline

\end{tabular}
\medskip

\noindent
{\bf Case vi} of a horizontal cusp $\horcusp$.
By Proposition~\ref{prop:CriterionVersality}
 we should prove that any germ $(x(r),z(r)):\R\to\R^2$ can be represented
 as a sum of vectors from the spaces $T_1^r$, $T_1^l$ and $N_1^p$, i.e.
 we solve the functional equations from the table
 $x(r)=2rf(r)+g(r^2,r^3)$ and $z(r)=-dr+cr^2+3r^2f(r)+h(r^3)$,
 which have one the of the possible solutions
$$\left\{ \begin{array}{l}
d=-\dot z(0),\quad h(r^3)=z(0), \quad
f(r)=\dfrac{ z(r)-\dot z(0)r-z(0) }{3r^2}+\dfrac{\dot x(0) }{2}-\dfrac{ \ddot z(0) }{6},\\
c=\dfrac{\ddot z(0)-3\dot x(0)}{2}, \;
g(r^2,r^3)=x(r)-\dot x(0)r-2\dfrac{ z(r)-\ddot z(0)r^2/2-\dot z(0)r-z(0) }{3r}
\end{array} \right.$$
Here $h$ has only the constant term and $g(r^2,r^3)$ has no linear term in $r$,
 all other powers have the form $2j+3k$ for some integers $j,k\geq 0$, e.g.
$$\mbox{  for a germ } (a_0+a_1r+a_2r^2+\dots,\quad b_0+b_1r+b_2r^2+\dots) \mbox{ one has }$$ 
$$f=a_1/2+\dots,\; g(x,z)=a_0+a_2x+\dots,\;
h=b_0,\; d=-b_1,\; c=(2b_2-3a_1)/2.$$

\noindent
{\bf Case vii} of a mixed tangency $\hortang$.
We prove that at each point $r_i$, $i=1,2$,
 any germ $(x_i,z_i):\R\to\R^2$ can be represented as 
 a sum of vectors from $T_i^r$, $T_i^l$, $N_i^p$,
 i.e. in terms of suitable $c,d$ and $f$, $g$, $h$.
Write down the equations from the table above.
$$\left\{\begin{array}{ll}
x_1(r)=f_1(r)+g(r,r^2), & z_1(r)=2rf_1(r)+h(r^2)-d,\\
x_2(r)=c+f_2(r)+g(r,-r^2), & z_2(r)=-2rf_2(r)+h(-r^2).
 \end{array} \right.
\leqno{(\hortang)}$$

For a function $f(r)$ denote its constant term simply by $f$.
The equations $z_1(r)=2rf_1(r)+h(r^2)-d$ and $z_2(r)=2rf_2(r)+h(-r^2)$ in degree~1 
 determine the constant terms $f_1,f_2$ of $f_1(r),f_2(r)$.
Then system $(\hortang)$ in degree~0 has a unique solution:
$$
\left\{\begin{array}{ll}
x_1=f_1+g, & z_1=h-d,\\
x_2=c+f_2+g, & z_2=h.
\end{array} \right.\;
\left\{\begin{array}{ll}
g=x_1-f_1, & h=z_2,\\
c=x_2-x_1+f_1-f_2, & d=z_2-z_1.
\end{array} \right.
\leqno{(\hortang_0)}$$
For a function $f(r)$ define its odd and even part as
$\odd f(r)=\dfrac{f(r)-f(-r)}{2}$,\\ $\even f(r)=\dfrac{f(r)+f(-r)}{2}$.
We look for solutions $g(x,z)=g_1(x)+g_2(z)$ 
 and $h(z)$ such that $g_2(-z)=-g_2(z)$, $h(z)=h(-z)$.
Split each equation of $(\hortang)$:
$$\left\{\begin{array}{l}
\odd x_1(r)=\odd f_1(r)+\odd g_1(r), \quad
\even z_1(r)=2r\odd f_1(r)+h(r^2)-d,\\
\even x_1(r)=\even f_1(r)+\even g_1(r)+g_2(r^2), \quad
\odd z_1(r)=2r\even f_1(r), \\
\odd x_2(r)=\odd f_2(r)+\odd g_1(r), \quad
\even z_2(r)=-2r\odd f_2(r)+h(r^2) \\
\even x_2(r)=c+\even f_2(r)+\even g_1(r)-g_2(r^2), \quad
\odd z_2(r)=-2r\even f_2(r).
\end{array} \right.$$
The resulting system has a solution below, where
 $\even z_1(r)-\even z_2(r)+d$ is divisible by $r$ due to~($\hortang_0$).
So the deformation is versal by Proposition~\ref{prop:CriterionVersality}.
$$\left\{\begin{array}{l}
\even f_1(r)=\odd z_1(r)/2r,\quad 
\even f_2(r)=-\odd z_2(r)/2r,\\
\odd f_1(r)=(\even z_1(r)-\even z_2(r)+d)/4r+(\odd x_1(r)-\odd x_2(r))/2,\\
\odd f_2(r)=(\even z_1(r)-\even z_2(r)+d)/4r+(\odd x_2(r)-\odd x_1(r))/2,\\
\even g_1(r)=(\even x_1(r)+\even x_2(r)-\odd z_1(r)/2r+\odd z_2(r)/2r-c)/2,\\
\odd g_1(r)=(+\even z_2(r)-\even z_1(r)-d)/4r+(\odd x_1(r)+\odd x_2(r))/2,\\
g_2(r^2)=(\even x_1(r)-\even x_2(r)-\odd z_1(r)/2r-\odd z_2(r)/2r+c)/2,\\
h(r^2)=(\even z_1(r)+\even z_2(r)+d)/2+r(\odd x_2(r)-\odd x_1(r)).
\end{array} \right.$$

\noindent
{\bf Case viii} of an extreme tangency $\maxtang$ is similar to Case vii.
\smallskip

\noindent
{\bf Case ix} of a horizontal triple point $\hortrip$.
The table above gives
$$\left\{\begin{array}{ll}
x_1(r)=f_1(r)+g(r,-r^2), & z_1(r)=-2rf_1(r)+h(-r^2),\\
x_2(r)=c+f_2(r)+g(r,r), & z_2(r)=f_2(r)+h(r), \\
x_3(r)=-d-f_3(r)+g(-r,r),& z_3(r)=f_3(r)+h(r).
 \end{array} \right.
\leqno{(\hortrip)}$$
The equation $z_1(r)=-2rf_1(r)+h(-r^2)$ in degree~1 determines
 the constant term $f_1$ of the function $f_1(r)$.
Then system $(\hortrip)$ in degree~0 has a unique solution.
$$\left\{\begin{array}{ll}
x_1=f_1+g, & z_1=h,\\
x_2=c+f_2+g, & z_2=f_2+h, \\
x_3=-d-f_3+g,& z_3=f_3+h.
\end{array} \right.\;
\left\{\begin{array}{ll}
g=x_1-f_1, & h=z_1,\\
f_2=z_2-z_1, & c=x_2-x_1+f_1+z_1-z_2, \\
f_3=z_3-z_1,& d=x_1-f_1-x_3+z_1-z_3.
\end{array} \right.$$

We look for $g(x,z)=g_1(x)+g_2(z)$.
Apply elementary operations to $(\hortrip)$
$$\left\{\begin{array}{l}
2rx_1(r)+z_1(r)=2rg_1(r)+2rg_2(-r^2)+h(-r^2),\\
x_2(r)-z_2(r)=c+g_1(r)+g_2(r)-h(r),\\
x_3(r)+z_3(r)=-d+g_1(-r)+g_2(r)+h(r).
 \end{array} \right.
\leqno{(\hortrip_1)}$$
The functions $f_1,f_2,f_3$ can be expressed 
 in terms of the solutions of $(\hortrip_1)$. 
Split the 1st equation of $(\hortrip_1)$ into the odd and even parts, 
 then apply operations to $(\hortrip_1)$:
$$2r\odd x_1(r)+\even z_1(r)=2r\odd g_1(r)+h(-r^2),$$
$$2r\even x_1(r)+\odd z_1(r)=2r\even g_1(r)+2rg_2(-r^2),\leqno{(\hortrip_2)}$$
$$x_3(r)+z_3(r)+x_2(r)-z_2(r)-2\even x_1(r)-\dfrac{\odd z_1(r)}{r}=
 c-d+2g_2(r)-2g_2(-r^2),$$
 which determines the coefficients of $g_2(r)=\sum_{i=0}^{\infty} e_i r^i$
 splitting into parts as follows.
Taking the odd part, we compute $e_i$ with all odd $i$,
 the consider terms with powers $4i$ and $4i+2$ separately,
 find all $e_{4i+2}$ and continue splitting into parts.
Having found $g_2(r)$, compute $\even g_1(r)$ from $(\hortrip_2)$
 and work out $h(r)$, $\odd g_1(r)$ from
$$\left\{\begin{array}{l}
x_2(r)-z_2(r)-x_3(r)-z_3(r)=c+d+2\odd g_1(r)-2h(r),\\
2r\odd x_1(r)+\even z_1(r)=2r\odd g_1(r)+h(-r^2),
 \end{array} \right.$$
 excluding $\odd g_1(r)$ and then splitting the result into parts as above.
\smallskip

\noindent
{\bf Case x} of another horizontal triple point $\inctrip$ is similar to Case~ix.
\qed
\end{sketch}


\subsection{Bifurcation diagrams of codimension~2 singularities}
\label{subs:BifurcationDiagrams}
\noindent
\smallskip

\noindent
The \emph{bifurcation} diagram of a codimension~2 singularity
 $\de$ from Definition~3.4 is formed by the pairs $(a,b)\in\R^2$
 from the versal deformation of $\de$ from Lemma~\ref{lem:VersalDeformations}.
We will describe curves representing codimension~1 subspaces $\Si_{\ga}$
 adjoined to $\Si_{\de}$ in the space $\sk$ of all links $K\subset V$.
\medskip

Oriented arcs in bifurcation diagrams of Figure~\ref{fig:BifurcationDiagrams} 
 are associated to canonical loops $\cl(K_{\pm\e})\subset\sk$, 
 where links $K_{\pm\e}$ are close to a given link $K_{0}$.
At the zero critical moment, the loop $\cl(K_0)$ defines
 an arc through the origin $\{a=b=0\}$.
These arcs are transversal to the codimension~1 subspace $\Si^{(1)}$ apart from the cases below.
In Figure~\ref{fig:BifurcationDiagrams}ix and \ref{fig:BifurcationDiagrams}x 
 the canonical loop $\cl(K_s)$ is \emph{parallel}
 to $\Si_{\maximin}$, $\Si_{\maximax}$, $\Si_{\cubicsn}$
 in the following sense: if $K\in\Si_{\maximin}$, then
 $\cl(K)\subset\Si_{\maximin}\cup\Si_{\hortang}$.
If $K\in\Si_{\cube}$, then
 $\cl(K)\subset\Si_{\cube}\cup\Si_{\horcusp}$.
Similarly, $K\in\Si_{\maximax}$ implies that
 $\cl(K)\subset\Si_{\maximax}\cup\Si_{\maxtang}$.
\medskip

\begin{lemma}
\label{lem:BifurcationDiagrams}
Figure~\ref{fig:BifurcationDiagrams} contains the bifurcation diagrams
 of the codimension~2 singularities $\de$ :
$\quadrup,\tangint,\cuspint,\cubicsn,\degcusp,
 \hortrip,\inctrip,\hortang,\horcusp,\maxtang$
 and shows how the canonical loops $\cl(K_{\pm\e})$
 intersect the adjoined codimension~1 subspaces $\Si_{\ga}$.
\end{lemma}
\begin{proof}
In Cases i--v below the canonical loops transversally intersects all 
 the singular subspaces since the tangents of intersecting arcs are not horizontal.
\smallskip

\begin{figure}
\includegraphics[scale=1.0]{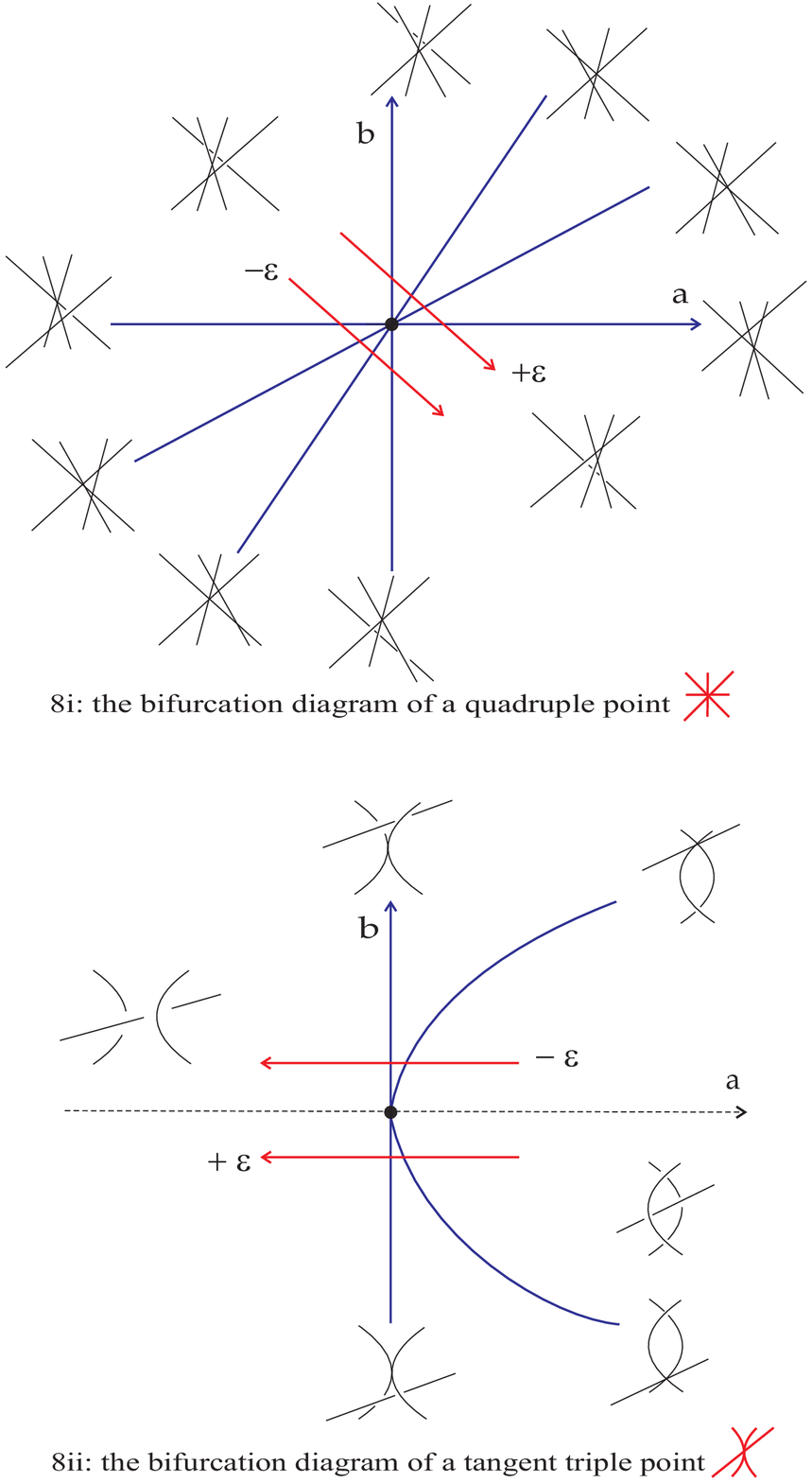}
\end{figure}

\begin{figure}
\includegraphics[scale=1.0]{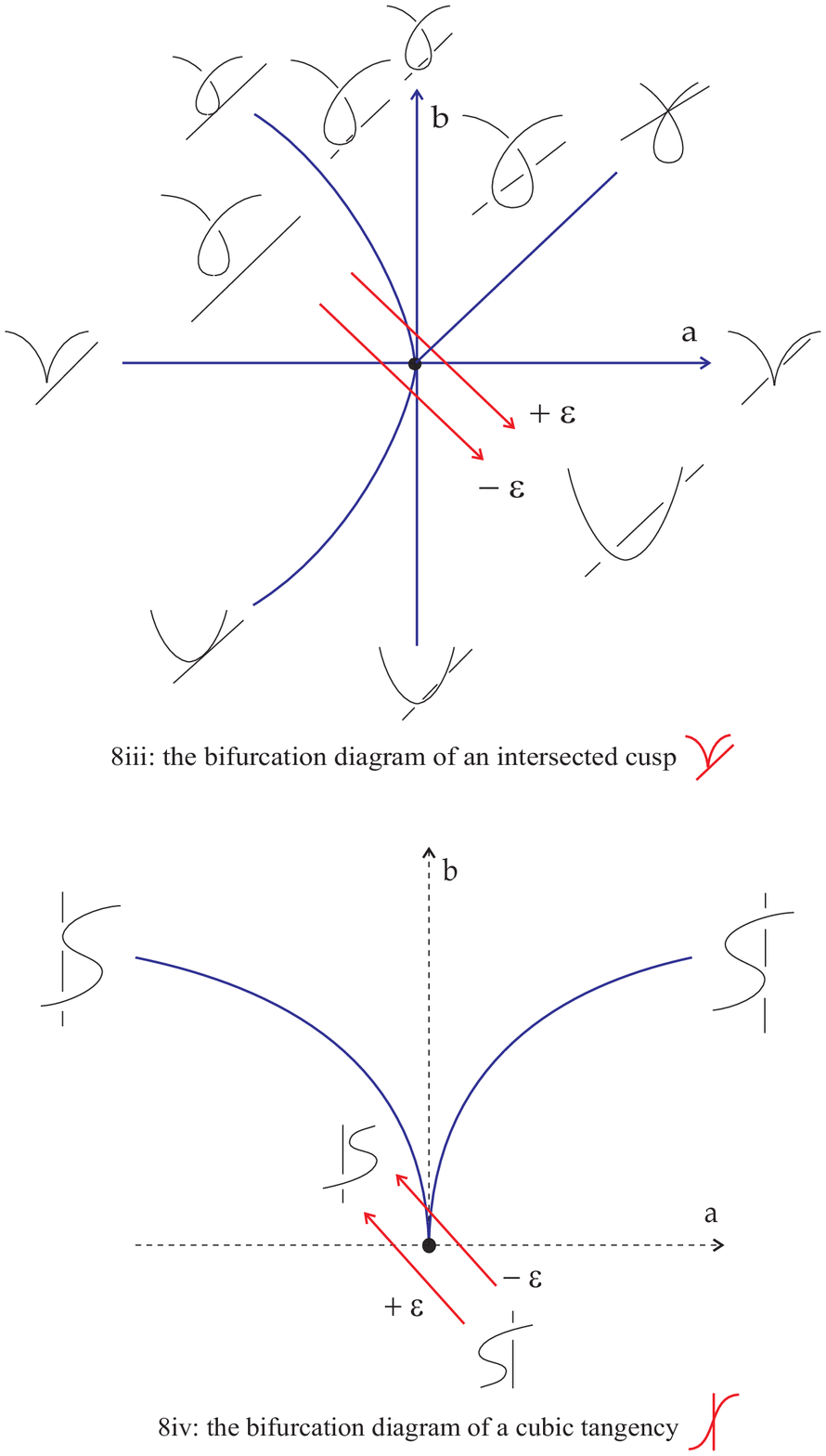}
\end{figure}

\begin{figure}
\includegraphics[scale=1.0]{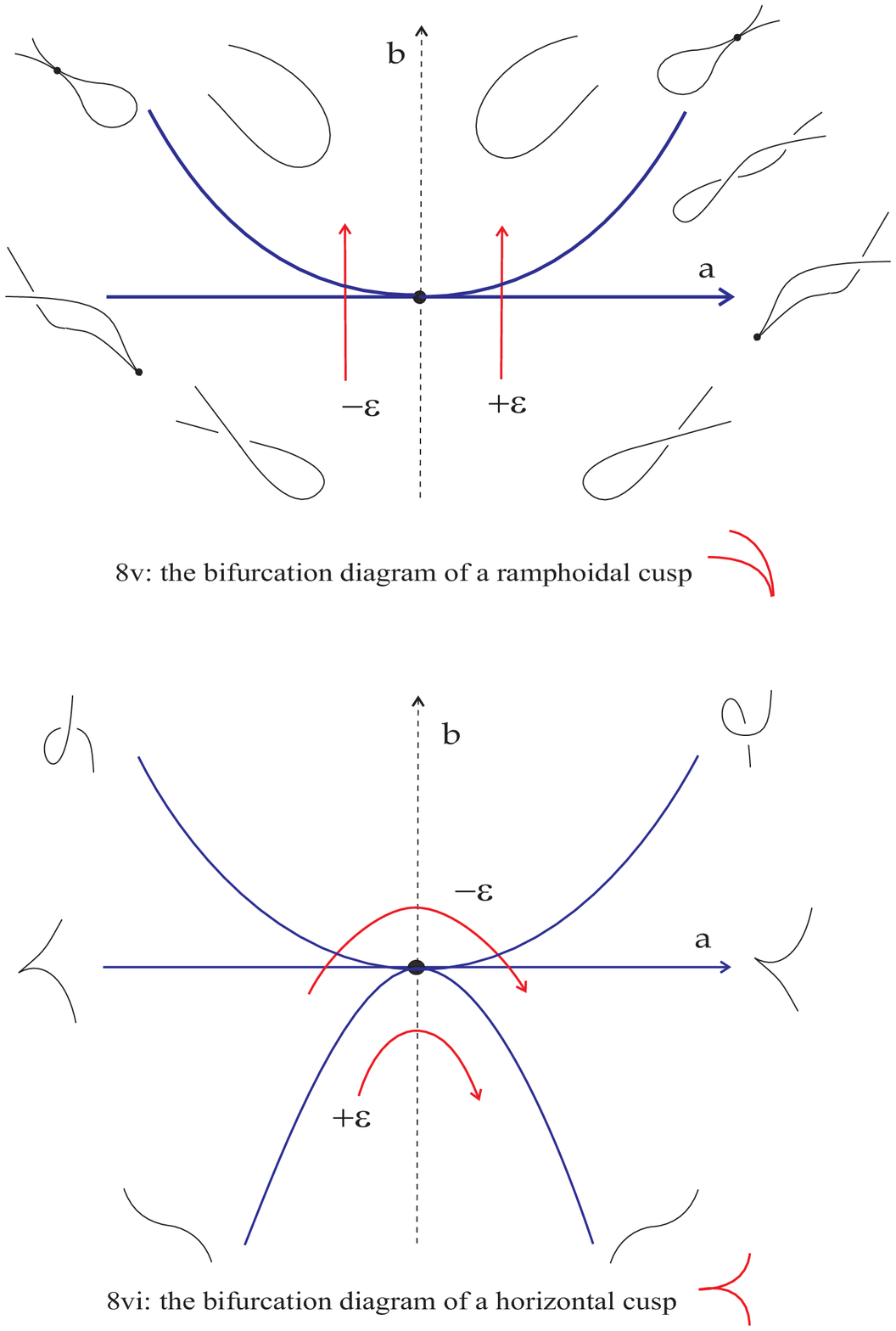}
\end{figure}

\begin{figure}
\includegraphics[scale=1.0]{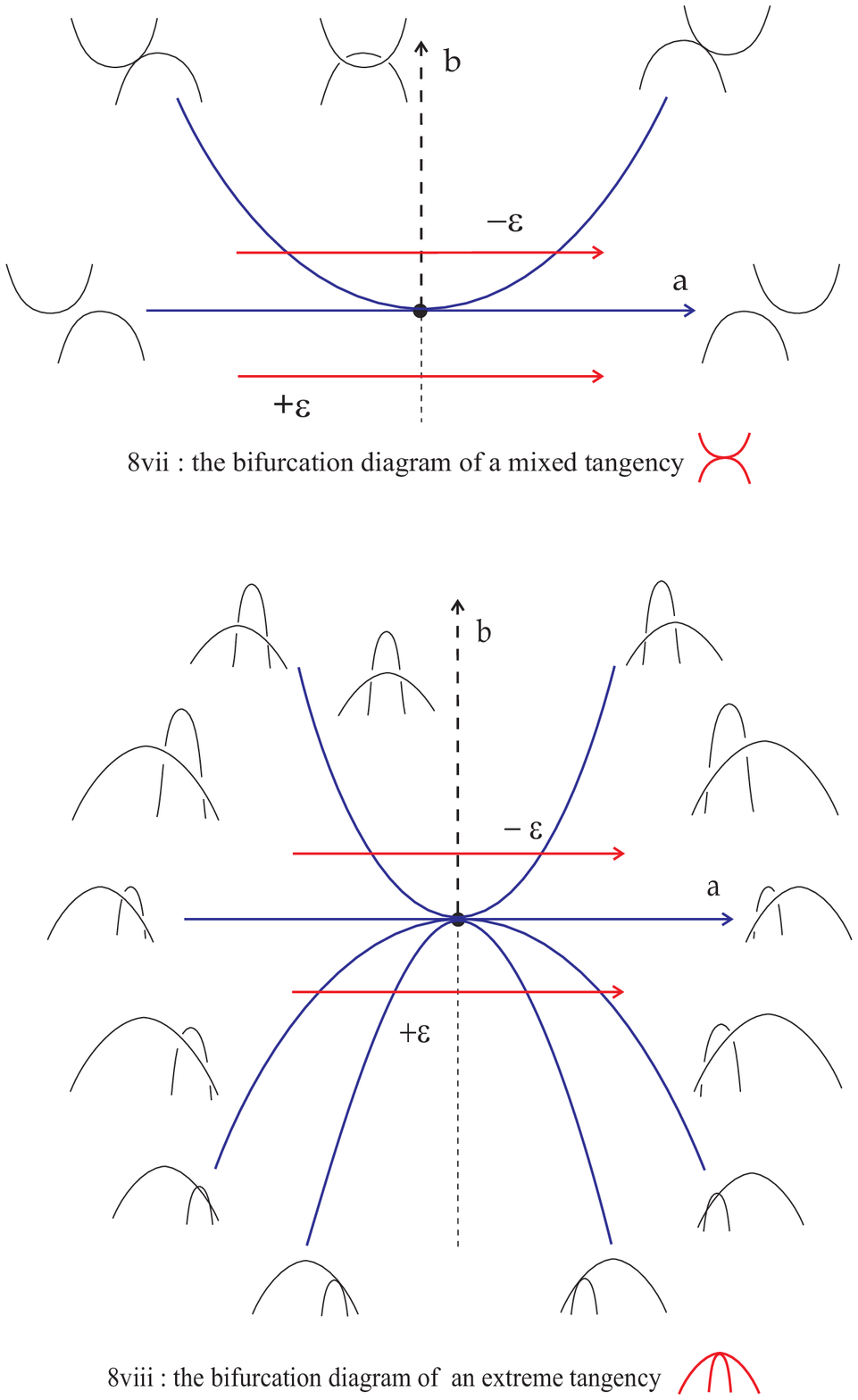}
\end{figure}

\begin{figure}
\includegraphics[scale=1.0]{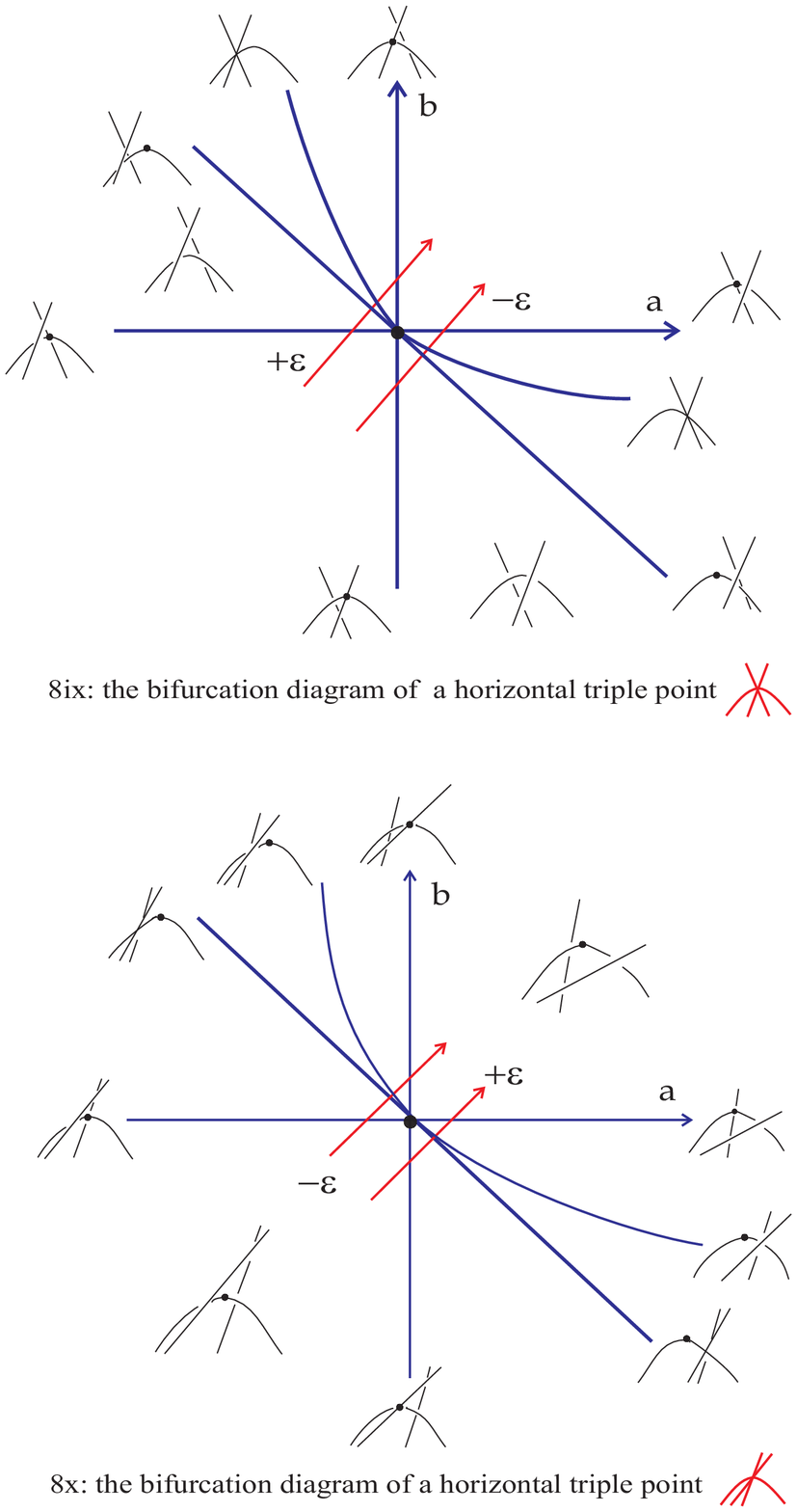}
\caption{Bifurcation diagrams of codimension~2 singularities}
\label{fig:BifurcationDiagrams}
\end{figure}

\noindent
{\bf Case i} of a quadruple point $\quadrup$.
There are 4 singular subspaces $\Si_{\trip}$ intersecting each other 
 transversally at the singular subspace $\Si_{\quadrup}$.
Using the normal form of $\quadrup$ from Lemma~\ref{lem:NormalForms}, 
 we show 4 subspaces in the bifurcation diagram of 
 Figure~\ref{fig:BifurcationDiagrams}i, namely
 $\{a=0\}$ (branches 1, 2, 4 intersect), $\{b=0\}$ (branches 1, 3, 4 intersect),
 $\{a=b\}$ (branches 1, 2, 3 intersect), $\{e(a+b)=b-a\}$ (branches 2, 3, 4 intersect). 
\smallskip

\noindent
{\bf Case ii} of a tangent triple point $\tantrip$.
The branches $\{x=z^2-a\}$, $\{x=0\}$ have a tangency if $a=0$.
The triple point appears when $z^2-a=0=z-b$, i.e. $a=b^2$.
The bifurcation diagram of Figure~\ref{fig:BifurcationDiagrams}ii
 has 1 parabola and 1 line touching each other.
\smallskip

\noindent
{\bf Case iii} of an intersected cusp $\cuspint$.
The branch $(r^3-br,r^2)$ has a self-intersection at $r=\pm\sqrt{b}$, 
 $b\geq 0$, which becomes an ordinary cusp if $b=0$.
The self-intersection is a triple point 
 when it is on the branch $(r-a,r)$, i.e. $a=b$.
Finally, we get a simple tangency of $(r^3-br,r^2)$ and $(r-a,r)$ 
 if $a=2r^3-r^2,b=3r^2-2r$ or $3a-2br=r^2$ has a double root, i.e. $b^2+3a=0$.
The bifurcation diagram of Figure~\ref{fig:BifurcationDiagrams}iii
  contains 1 parabola, 1 line and 1 ray meeting at 0.
\smallskip

\noindent
{\bf Case iv} of a cubic tangency $\cubicsn$.
The branch $(r^3-3br+a,r)$ has extrema of the $x$-coordinate
 at $r=\pm\sqrt{b}$, which lie on $(0,r)$ 
 if $r^3-3br+a=0$, i.e. $a^2=4b^3$.
The only subspace $\Si_{\tang}$ is adjoined to $\Si_{\cubicsn}$
 in the bifurcation diagram of Figure~\ref{fig:BifurcationDiagrams}iv.
\smallskip

\noindent
{\bf Case v} of a ramphoidal cusp $\degcusp$.
The curve $(r^5+ar^3+br,r^2)$ has an ordinary cusp 
 when $\dot x=\dot z=0$, ie $r=0$ and $b=0$, and 
 a self-tangency when $5r^4+3ar^2+b=0$ has two double roots, i.e. $9a^2=20b$.
The bifurcation diagram of Figure~\ref{fig:BifurcationDiagrams}v
 contains 1 parabola and 1 line touching each other at 0.
\smallskip

\noindent
{\bf Case vi} of a horizontal cusp $\horcusp$.
The curve $(r^2,r^3+ar^2-br)$ has a crossing at $\pm r$, hence $r^3=br$ and $r=\pm\sqrt{b}$, $b>0$.
This crossing is critical, ie $\dot z=3r^2+2ar-b=0$, if $b=a^2$.
The critical point becomes degenerate, ie $\ddot z=6r+2a=0$, if $b=-a^2/3$.
The subspace $\Si_{\cusp}$ of ordinary cusps, 
 where $\dot x=\dot z=0$, is represented by $\{b=0\}$.
The bifurcation diagram of Figure~\ref{fig:BifurcationDiagrams}vi 
 shows 2 parabolas, 1 line and 1 ray meeting at 0.
The arc associated to a canonical loop moves in the vertical direction
 and remains parallel to the parabola $\{b=-a^2/3\}$ representing the subspace $\Si_{\cube}$.
\smallskip

\noindent
{\bf Case vii} of a mixed tangency $\hortang$.
The branch $(r,r^2-b)$ touches $(r+a,-r^2)$ if 
 $r^2-b=-(r-a)^2$ has a double root, i.e. $a^2=2b$.
Both curves have extrema in the same horizontal line when $b=0$.
The bifurcation diagram of Figure~\ref{fig:BifurcationDiagrams}vii has 
 1 parabola and 1 line touching each other at 0.
\smallskip

\noindent
{\bf Case viii} of an extreme tangency $\maxtang$.
The branch $(r,-2r^2+b)$ touches $(r+a,-r^2)$ if 
b$-2r^2+b=-(r-a)^2$ has a double root, i.e. $2a^2+b=0$.
Both branches have extrema in the same horizontal line when $b=0$.
The branch $(r,-2r^2+b)$ passes through an extremum of $(r+a,-r^2)$ at $r=0$ if $b=2a^2$.
The branch $(r+a,-r^2)$ passes through an extremum of $(r,-2r^2+b)$ at $r=0$ if $b=-a^2$.
The bifurcation diagram of Figure~\ref{fig:BifurcationDiagrams}vii has 
 3 parabolas and 1 line touching each other at 0.
\smallskip

\noindent
{\bf Case ix} of a horizontal triple point $\hortrip$.
The branches $(r+a,r)$ and $(-r-b,r)$ pass through the extremum of $(r,-r^2)$ at $r=0$
 when $a=0$ and $b=0$, respectively.
The crossing of $(r+a,r)$ and $(-r-b,r)$ at $r=-(a+b)/2$ lies in the same horizontal line
 with the extremum of $(r,-r^2)$ at $r=0$ if $a+b=0$.
The branches $(r,-r^2)$, $(r+a,r)$ and $(-r-b,r)$ have a triple point if
 $r=-r^2+a=r^2-b$ or $(a-b)^2=2(a+b)$, which is a parabola
 in the bifurcation diagram of Figure~\ref{fig:BifurcationDiagrams}ix.
The arc associated to a canonical loop is transversal to the subspaces,
 because only one tangent remains horizontal under the rotation.
\smallskip

\noindent
{\bf Case x} of another horizontal triple point $\inctrip$ is similar to Case ix.
\end{proof}


\section{The diagram surface of a link}
\label{sect:DiagramSurfaceOfLink}

In this section the classification problem of
 generic links $K\subset V$ reduces to
 their diagram surfaces $\ds(K)$
 in the thickened torus $\T=\axz\times S_t^1$,
 $\axz=[-1,1]_x\times S_z^1$.
\smallskip


\subsection{The diagram surface of a link and generic surfaces}
\label{subs:GenericSurfaces}
\noindent
\smallskip

Briefly the diagram surface of
 a loop $\{K_t\}$ of links is the 1-parameter family of
 the diagrams $\pr_{xz}(K_t)\subset\axz\times\{t\}$.
This family can be considered as the union of link diagrams,
 i.e. as a 2-dimensional surface in the thickened torus $\T=\axz\times S_t^1$.
\smallskip

\begin{definition}
\label{def:DiagramSurfaceOfLink}
Let $\{K_t\}\subset\sk$ be a loop of links.
The \emph{diagram surface} $\ds(\{K_t\})\subset\axz\times S_t^1$
 is formed by the diagrams
 $\pr_{xz}(K_t)\subset\axz\times\{t\}$, $t\in S_t^1$.
If $K_t$ are knots, $\ds(\{K_t\})$ is the torus $S^1\times S_t^1$ mapped to
 the thickened torus $\T=\axz\times S_t^1$.
The \emph{diagram surface} $\ds(K)$ of
 an oriented link $K\subset V$ consists of the diagrams
 $\pr_{xz}(\rot_t(K))\subset\axz\times\{t\}$
 and is oriented by the orientations of $K$ and $S_t^1$.
\ed
\end{definition}

\begin{figure}[!h]
\includegraphics[scale=1.0]{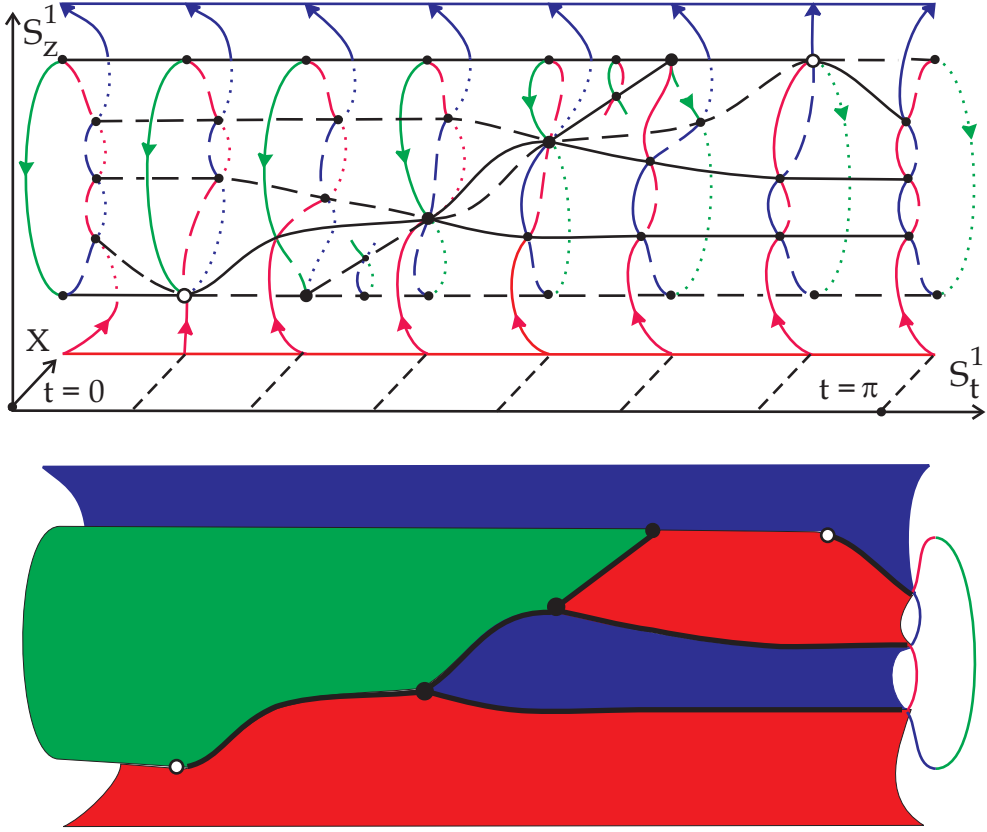}
\caption{Half the diagram surface of a smoothed trefoil 
 from Figure~\ref{fig:RotatedTrefoils}.}
\label{fig:DiagramSurfaceTrefoil}
\end{figure}

Figure~\ref{fig:DiagramSurfaceTrefoil} shows vertical sections of $\ds(K)$
 for a smoothed trefoil $K$ from Figure~\ref{fig:RotatedTrefoils}, $t\in[0,\pi]$.
Each section is the diagram of a rotated trefoil $\rot_t(K)$
 for some $t\in S_t^1$.
Local extrema of $\rot_t(K)$ form horizontal circles parallel
 to $S_t^1$.
Several arcs in Figure~\ref{fig:DiagramSurfaceTrefoil} are dashed or dotted, because
 they are invisible in the $x$-direction.
\smallskip

By Definition~\ref{def:GenericLoop} the shift $t\mapsto t+\pi$ maps 
 the surface $\ds(\{K_t\})$ to its image under the symmetry in $S_z^1$.
Actually the link $K_{t+\pi}$ is obtained from $K_t$
 by the symmetry $\rot_{\pi}$, ie
 the diagrams $\pr_{xz}(K_{t+\pi})$ and
 $\pr_{xz}(K_t)$ are symmetric for all $t\in S_t^1$.
For a generic loop $\{K_t\}$, the vertical sections of $\ds(K)$
 are the diagrams $\pr_{xz}(K_t)$
 and allow the codimension~1 singularities
 $\trip,\tang,\cusp,\crit$ only.
It follows from the fact that any critical point
 of $\pr_z:K_t\to S_z^1$ remains critical under $\rot_t$.
\smallskip

For any $t\in S_t^1$, the points from $K_t\cap(\dxy\times\{z=\pm 1\})$ 
 and the critical points of $\pr_z:K_t\to S_z^1$ 
 divide the $i$-th component of $K_t$ into arcs $A_{t,i,q}$, $q=1,\dots,n_i$.
The total number of these arcs does not depend on $t$
 since any critical point $a_t\in K_t$ of $\pr_z$
 remains critical while $t$ varies.
The union $\cup a_t$ of the extrema of $\pr_z:K_t\to S_z^1$
 for all $t\in S_t^1$ splits into
 \emph{critical circles} $C_i$ of $\ds(\{K_t\})$.
The union $B_{i,q}=\cup A_{t,i,q}$ over all $t\in S_t^1$
 is called a \emph{trace band} of $\ds(\{K_t\})$.
The 3 trace bands in the bottom picture of Figure~\ref{fig:DiagramSurfaceTrefoil}
 have different colours.
The arcs $A_{t,i,q}$ are monotonic with respect
 to $\pr_{zt}:K_t\to S_z^1\times\{t\}$.
Then the trace bands project
 1-1 under $\pr_{zt}:\ds(\{K_t\})\to S_z^1\times S_t^1$.
Successive bands $B_{i,q}$, $B_{i+1,q}$
 meet at a critical circle.
\smallskip

The \emph{singular} points of $\ds(\{K_t\})$
 are crossings and codimension~1 singularities of
 the diagrams $\pr_{xz}(K_t)$ over all $t\in S_t^1$.
\emph{A trace arc} is an intersection of the interiors of
 2 trace bands in $\ds(\{K_t\})$.
The triple points, tangent points, cusps
 and critical crossings of link diagrams
 $\pr_{xz}(K_t)$ are called \emph{triple} vertices,
 \emph{tangent} vertices, \emph{hanging} vertices and
 \emph{critical} vertices of $\ds(\{K_t\})$, respectively.
So a trace arc may contain several vertices of 
 $\ds(\{K_t\})$ in the usual sense.
\smallskip

Take a singular point $p\in\ds(K)$ that is not a vertex
 and does not belong to a critical circle of $\ds(K)$.
Then $p$ is a double crossing of two arcs
 $A_{t,i,q}$ and $A_{t,j,s}$ in a diagram $\pr_{xz}(K_t)$.
If the arc $A_{t,i,q}$ passes over (respectively, under)
 $A_{t,j,s}$ then associate to $p$ the \emph{label} $(q_i s_j)$
 (respectively, the \emph{reversed} label $(s_j q_i)$).
If $K_t$ is a knot 
 then we miss the indices $i,j=1$ as in Figure~\ref{fig:TraceGraphTrefoil}.
\smallskip

Trace arcs of $\ds(\{K_t\})$ end at hanging vertices,
 meet each other at critical vertices
 and intersect at triple vertices.
Each trace arc of $\ds(K)$ is the evolution trace
 of a double crossing in $\rxz\times S_t^1$ while $t$ varies.
The label of a point $p$ does not change when $p$
 passes through tangent vertices and triple vertices.
\smallskip

The diagram surface can be defined for
 any loop of links and can be extremely complicated.
The surfaces corresponding to generic loops
 are simple and play the role
 of general link diagrams in dimension~3.
As in the case of links, we define a generic surface
 associated to a generic loop.
A generic surface will be an immersed surface
 with all combinatorial features
 of diagram surfaces of generic loops.
For any generic surface, a corresponding generic loop
 is constructed in Lemma~\ref{lem:ReconstructGenericLoop}.
\smallskip

\begin{definition}
\label{def:GenericSurface}
Decompose $S_i^1$ into arcs $A_{i,1},\dots,A_{i,n_i}$.
Introduce the \emph{trace bands} $B_{i,q}=A_{i,q}\times S_t^1$,
 $q=1,\dots,n_i$.
A \emph{generic} surface $S$ is the image of a smooth map
 $h:(\sqcup_{i=1}^m S_i^1)\times S_t^1=\cup_{i=1}^m (\cup_{q=1}^{n_i} B_{iq})\to\rxz\times S_t^1$
 such that Conditions (i)--(v) hold
\medskip

\noindent
{\bf (i) Conditions} on \emph{symmetry} and \emph{trace bands}.
\medskip

$\bu$
 under $t\mapsto t+\pi$ the surface $S$ maps to
 its image under the symmetry in $S_z^1$;
\smallskip

$\bu$
 each trace band $B_{i,q}\subset S$ projects one-to-one
 under $\pr_{zt}:S\to S_z^1\times S_t^1$.
\medskip

\noindent
The surface $S$ should be simple enough.
More formally we require the following.
\bigskip

\noindent
{\bf (ii) Conditions} on \emph{sections}
  $D_t=S\cap(\axz\times\{t\})$, $t\in S_t^1$.
\medskip

\noindent
There are finitely many critical moments
 $t_1,\dots,t_l\in S_t^1$ such that
\smallskip

$\bu$
for all $t\notin\{t_1,\dots,t_l\}$,
 the sections $\{D_t\}$ are general diagrams;
\smallskip

$\bu$
for each $t=t_1,\dots,t_l$,
 the section $D_t$ has one of the singularities
 $\trip,\tang,\cusp,\crit$;
\smallskip

$\bu$
while $t$ passes a critical moment,
 $D_t$ changes by a move I--IV in Figure~\ref{fig:ReidemeisterMoves}.
\bigskip

\noindent
Conditions (ii) on sections imply
 some restrictions on trace bands.
These requirements can be stated independently
 to define trace arcs and critical circles.
\medskip

\noindent
{\bf (iii) Conditions} on \emph{trace arcs} and
 \emph{critical} circles:
\medskip

$\bu$
 a \emph{trace arc} is an intersection of
 the interiors of 2 trace bands $B_{i,q}$ and $B_{j,s}$;
\smallskip

$\bu$
 a \emph{critical} circle $C_{i,q}$ is the common boundary
 of successive bands $B_{i,q}$, $B_{i,q+1}$;
\bigskip

\noindent
The arcs defined above allow us to introduce
 vertices of a generic surface $S$.
\medskip

\noindent
{\bf (iv) Conditions} on \emph{vertices}:
\medskip

$\bu$
 a \emph{triple} vertex is a transversal intersection of
 3 trace bands $B_{iq},B_{js},B_{kr}$;
\smallskip

$\bu$
 a \emph{hanging} vertex of $S$ is the endpoint
 of a trace arc in $B_{i,q}\cap B_{i,q+1}$;
\smallskip

$\bu$
 a \emph{critical} vertex is the intersection
 of a critical circle $C_{i,q}$ and $B_{j,s}\not\supset C_{i,q}$;
\smallskip

$\bu$
 a \emph{tangent} vertex is a critical point of
 $\pr_t$ on the interior of a trace arc;
\smallskip

$\bu$
 all the \emph{vertices} are distinct
 and map on different points under $\pr_t:S\to S_t^1$.
\bigskip

\noindent
Finally fix labels $(i,q)$ and $(j,s)$.
Take a trace arc from the intersection 
 $B_{i,q}\cap B_{j,s}$ of interiors of 2 trace bands.
Endow the chosen arc with a \emph{label}:
 either $(q_i s_j)$ or $(s_j q_i)$ in such a way that
 the following restrictions apply.
\medskip

\noindent
{\bf (v) Conditions} on \emph{labels}:
\medskip

$\bu$
 under the time shift $t\mapsto t+\pi$,
 each label reverses: $(q_i s_j)\mapsto(s_j q_i)$;
\smallskip

$\bu$
 trace arcs intersecting at a triple vertex
 are endowed with $(q_i s_j)$, $(s_j r_k)$, $(q_i r_k)$;
\smallskip

$\bu$
 a hanging vertex is endowed with the label
 of the trace arc containing it;
\smallskip

$\bu$
 each circle $C_{i,q}$ has 2 hanging vertices
 endowed with $((q+1)_i,q_i)$, $(q_i,(q+1)_i)$;
\smallskip

$\bu$
 if a trace band $B_{j,s}$ intersects a critical circle $C_{i,q}$
 in a vertex $c$ then the label\\
 \hspace*{6mm}
 at $c$ transforms as follows:
 $(q_i s_j)\lra((q+1)_i,s_j)$ or $(s_j q_i)\lra(s_j,(q+1)_i)$.
\ed
\end{definition}
\smallskip

To get the following result compare
 Definitions~\ref{def:GenericLoop}, \ref{def:DiagramSurfaceOfLink}
 with Definition~\ref{def:GenericSurface}.
\smallskip

\begin{lemma}
\label{lem:GenericLoopsSurfaces}
For any generic loop $L$ of links,
 the diagram surface $\ds(L)$ is a generic surface
 in the sense of Definition~\ref{def:GenericSurface}.
\qed
\end{lemma}
\smallskip


\subsection{Three-dimensional moves on generic surfaces}
\noindent
\smallskip

\begin{figure}
\includegraphics[scale=1.0]{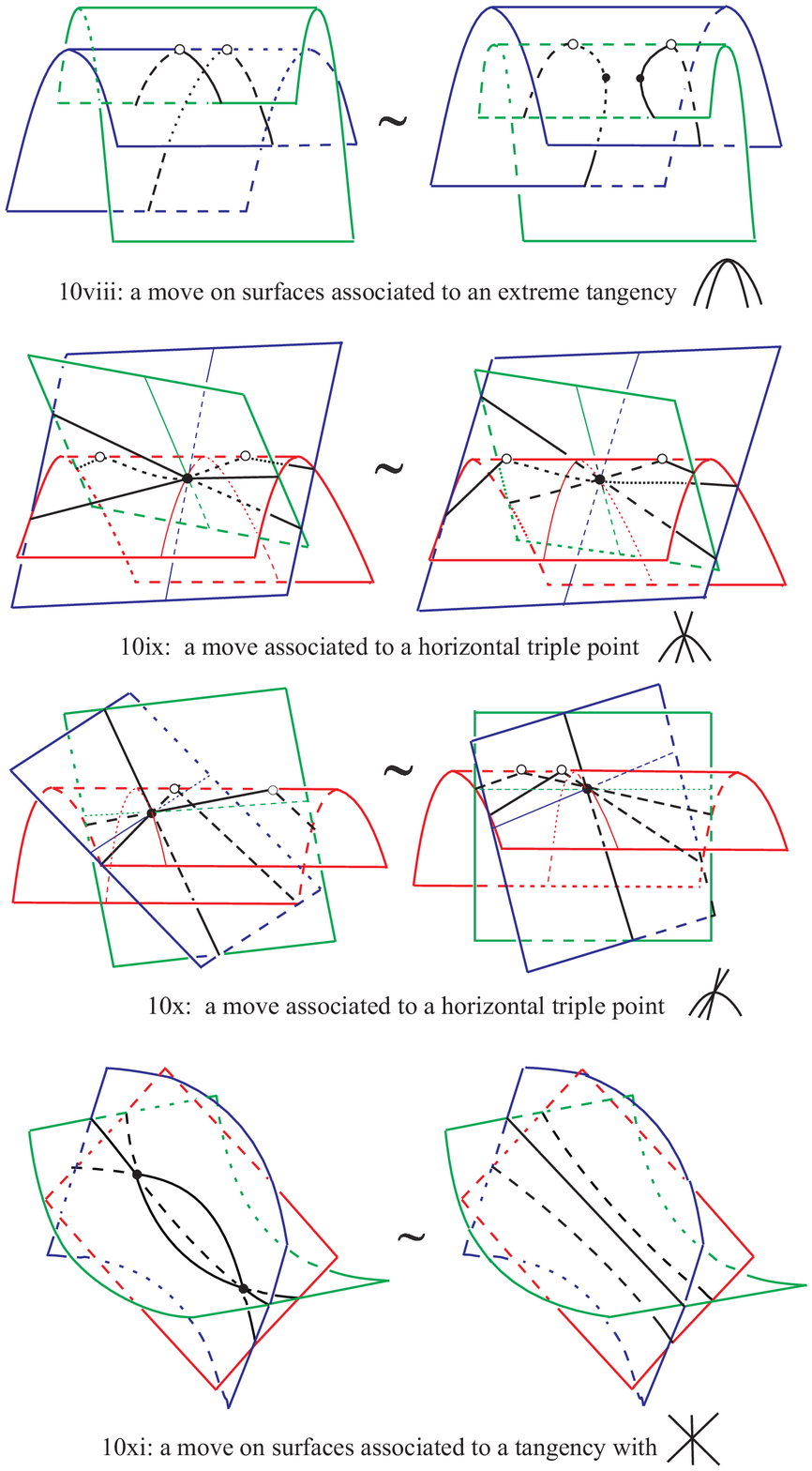}
\caption{Three-dimensional moves on diagram surfaces}
\label{fig:MovesDiagramSurfaces}
\end{figure}

\begin{definition}
\label{def:EquivalenceSurfaces}
A smooth family of surfaces $\{S_r\subset\axz\times S_t^1\}$,
 $r\in[0,1]$, is an \emph{equivalence} if
 there are finitely many critical moments
 $r_1,\dots,r_k\in[0,1]$ such that
\smallskip

$\bu$
 for all non-critical moments $r\notin\{r_1,\dots,r_k\}$,
 the surfaces $S_r$ are generic;
\smallskip

$\bu$
 if $r$ passes through a critical moment,
 $S_r$ changes by a move in Figure~\ref{fig:MovesDiagramSurfaces}.
\ed
\end{definition}

\begin{figure}
\includegraphics[scale=1.0]{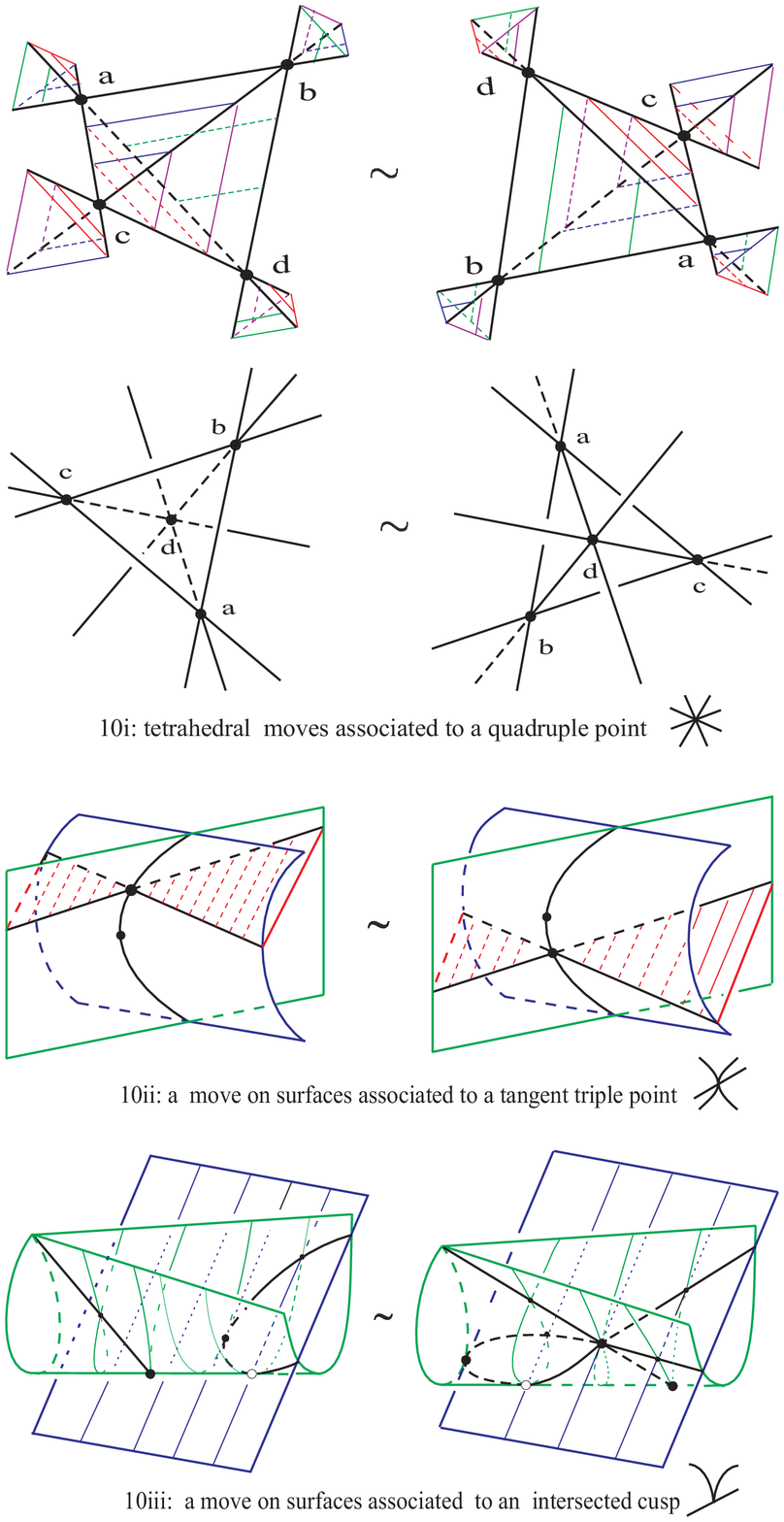}
\end{figure}

\begin{figure}
\includegraphics[scale=1.0]{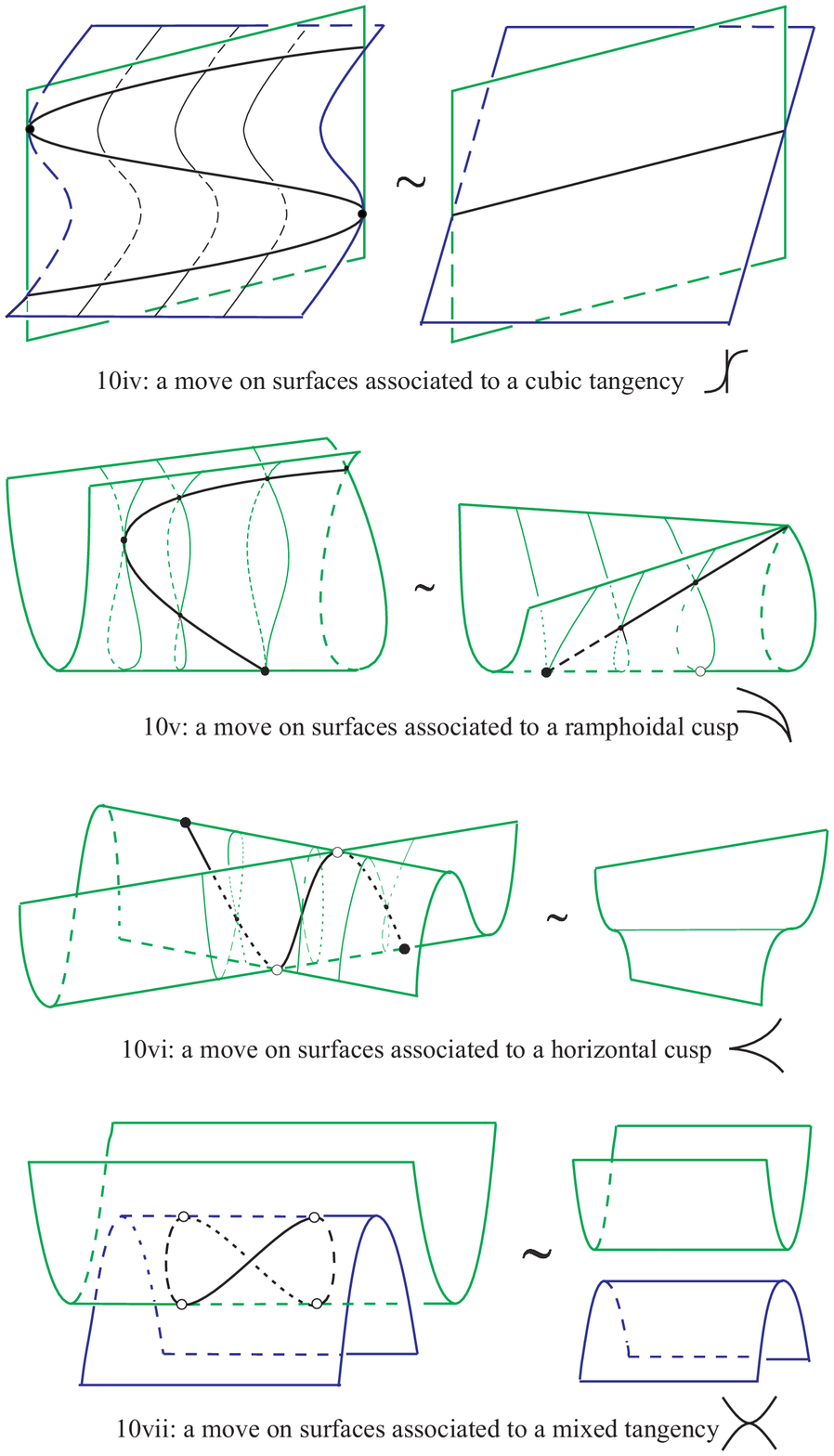}
\end{figure}

Each move in Figure~\ref{fig:MovesDiagramSurfaces} denotes 2 symmetric moves since
 the surfaces $S_r$ are symmetric in $S_z^1$ under $t\mapsto t+\pi$.
The following claim will be proved using bifurcation diagrams
 of codimension~2 singularities of link diagrams, see Lemma~\ref{lem:BifurcationDiagrams}.
\smallskip

\begin{lemma}
\label{lem:MovesDiagramSurfaces}
{\bf (a)}
Suppose that a family of loops $\{L_s\}$, $s\in[-1,1]$,
 in the space $\sk$ of all links $K\subset V$ transversally intersects
 the subspace $\Si^{(2)}$ at $s=0$.
Then the diagram surface $\ds(L_s)$ changes near
 $0$ by a move in Figure~\ref{fig:MovesDiagramSurfaces}i--x.
\medskip

\noindent
{\bf (b)}
If a family of loops $\{L_s\}$, $s\in[-1,1]$, in the space $\sk$ 
 has a simple tangency with $\Si_{\trip}$ at $s=0$, then
 $\ds(L_s)$ changes near $0$ by the move in 
 Figure~\ref{fig:MovesDiagramSurfaces}xi.
\end{lemma}
\begin{sketch}
The pictures in Figure~\ref{fig:MovesDiagramSurfaces} are obtained from
 the corresponding pictures in Figure~\ref{fig:BifurcationDiagrams}.
For instance, in Figure~\ref{fig:BifurcationDiagrams}iii 
 the canonical loop $\cl(K_{-\e})$
 meets 3 subspaces $\Si_{\cusp},\Si_{\tang},\Si_{\crit}$.
Therefore the surface $\ds(K_{-\e})$ has three distinguished
 points: a hanging vertex, a tangent vertex and a critical one
 as in Figure~\ref{fig:MovesDiagramSurfaces}iii.
Right after the move when all three points pass through each
 other, the surface  $\ds(K_{+\e})$ has 4 interesting points:
 three have the previous types, the new one is a triple vertex.
This situation agrees with 4 intersections of $\cl(K_{+\e})$
 with codimension~1 subspaces in Figure~\ref{fig:BifurcationDiagrams}iii.
The remaining cases are absolutely analogous.
\qed
\end{sketch}
\medskip

We produced Figure~\ref{fig:MovesDiagramSurfaces} first using
 our geometric intuition and then justified the moves
 applying the singularity theory in section~\ref{sect:ThroughCodim2Singularities}.
Since the family of sections in a generic surface
 is a general equivalence of diagrams then
 Lemma~\ref{lem:ReconstructGenericLoop} follows.
\smallskip

\begin{lemma}
\label{lem:ReconstructGenericLoop}
{\bf (a)}
For any generic surface $S$,
 there is a generic loop $L$ of links
 such that the diagram surface $\ds(L)$
 coincides with $S$.
\medskip

\noindent
{\bf (b)}
For any equivalence of surfaces $\{S_r\subset\axz\times S_t^1\}$,
 there is a generic homotopy of loops $\{L_r\}$
 such that $\ds(L_r)=S_r$, $r\in[0,1]$.
\qed
\end{lemma}
\smallskip

Lemma~\ref{lem:MovesDiagramSurfaces} and Definition~\ref{def:GenericHomotopy} 
 of a generic homotopy imply Lemma~\ref{lem:GenericHomotopyEquivalence}.
\smallskip

\begin{lemma}
\label{lem:GenericHomotopyEquivalence}
Any generic homotopy of loops $\{L_s\}$, $s\in[0,1]$
 in the space $\sk$ provides
 an equivalence $\{\ds(L_s)\}$ of diagram surfaces.
\qed
\end{lemma}
\smallskip

\begin{lemma}
\label{lem:EquivalenceGenericHomotopy}
Let $L_0,L_1$ be generic loops of links.
If $\ds(L_0)$ and $\ds(L_1)$ are equivalent in
 the sense of Definition~\ref{def:EquivalenceSurfaces}, then
 $L_0$ and $L_1$ are generically homotopic.
\end{lemma}
\begin{proof}
Any equivalence of diagram surfaces gives
 rise to a smooth family of loops $\{L_r\}$ by Lemma~\ref{lem:ReconstructGenericLoop}b.
The constructed family $\{L_r\}$ is a generic homotopy
 since all moves in Figure~\ref{fig:MovesDiagramSurfaces} correspond
 to singularities in the sense of Definition~\ref{def:Codim2Singularities}.
\end{proof}
\smallskip

By Lemmas~\ref{lem:GenericHomotopyEquivalence} and \ref{lem:EquivalenceGenericHomotopy} 
 the classification of generic links reduces to the equivalence problem for
 their diagram surfaces.
\smallskip

\begin{proposition}
\label{prop:ReductionToDiagramSurfaces}
Generic links $K_0,K_1$ are generically equivalent
 if and only if the diagram surfaces $\ds(K_0),\ds(K_1)$
 are equivalent in the sense of Definition~\ref{def:EquivalenceSurfaces}.
\qed
\end{proposition}
\smallskip

The isotopy class of a link can be easily
 reconstructed from its plane diagram,
 hence from its diagram surface with labels.
Formally, one has the following.
\smallskip

\begin{lemma}
\label{lem:ReconstructLinkSurface}
Suppose that the diagram surface $\ds(K)$ of
 a generic link $K$ is given, but $K$ is unknown.
Then one can reconstruct the isotopy class of $K\subset V$.
\qed
\end{lemma}
\smallskip


\section{The trace graph of a link as a link invariant}
\label{sect:TraceGraphLink}


\subsection{The trace graph of a link and generic trace graphs}
\noindent
\smallskip

Here the classification of links $K\subset V$
 will be reduced to their trace graphs.
\smallskip

\begin{definition}
\label{def:TraceGraphSurface}
Let $S\subset\rxz\times S_t^1$ be the diagram surface
 of a loop of links.
The \emph{trace graph} $\tg(S)$ is
 the self-intersection of $S$, ie
 a finite graph embedded into $\rxz\times S_t^1$.
The \emph{trace graph} $\tg(K)$ of a link $K$ is
 the trace graph of its diagram surface $\ds(K)$.
The trace arcs of $\ds(K)$ are called
 \emph{trace arcs} of $\tg(K)$.
The trace graph inherits the vertices and labels from $\ds(K)$.
\ed
\end{definition}
\smallskip

\begin{definition}
\label{def:GenericTraceGraph}
A finite graph $G\subset\rxz\times S_t^1$ is
 \emph{generic} if Conditions (i)--(ii) hold.
\medskip

\noindent
{\bf (i) Conditions} on \emph{trace arcs} and \emph{vertices}.
\smallskip

$\bu$
the graph $G$ consists of finitely many \emph{trace arcs},
 which are monotonic ars\\
 \hspace*{5mm}
 with respect to the orthogonal projection $\pr_z:G\to S_z^1$;
\smallskip

$\bu$
 any endpoint of a trace arc of $G$ has either\\
 \hspace*{5mm}
 degree~1 (a \emph{hanging} vertex $\hangver$) or
 degree~2 (a \emph{critical} vertex $\critver$);
\smallskip

$\bu$
 the critical vertices of $G$ coincide with the critical points of
 $\pr_z:G\to S_z^1$;
\smallskip

$\bu$
 trace arcs of $G$ intersect transversally at
 \emph{triple vertices} ($\tripver$);
\smallskip

$\bu$
 the critical points of $\pr_z:G\to S_t^1$
 are called \emph{tangent vertices} ($\tangver$).
\bigskip

\noindent
{\bf (ii) Conditions} on \emph{labels}.
\smallskip

$\bu$
 each trace arc of $G$ is labelled with a \emph{label} $(q_i s_j)$
 as in Definition~4.2;
\smallskip

$\bu$
under $t\mapsto t+\pi$ the graph $G$ maps
 to its image under the symmetry in $S_z^1$;
\smallskip

$\bu$
 under the time shift $t\mapsto t+\pi$
 each label $(q_i s_j)$ reverses to $(s_j q_i)$;
\smallskip

$\bu$
 every triple vertex $v\in G$ is labelled with
 \emph{a triplet} $(q_i s_j)$, $(s_j r_k)$, $(q_i r_k)$\\
 \hspace*{5mm}
 consisting of the labels
 associated to the trace arcs passing through $v$;
\smallskip

$\bu$
 each hanging vertex is labelled with
 the label of the corresponding trace arc;
\smallskip

$\bu$
 for any $i$ and $q=1,\dots,n_i$, there are
 exactly two hanging vertices of $G$\\
 \hspace*{5mm}
 labelled with $((q+1)_i,q_i)$ and $(q_i,(q+1)_i)$, respectively;
\smallskip

$\bu$
 at every critical vertex of $G$ the labels of trace arcs
 may transform as follows:\\
 \hspace*{5mm}
 either $(q_i s_j)\lra(q_i,(s\pm 1)_j)$ or $(q_i s_j)\lra((q\pm 1)_i,s_j)$.
\ed
\end{definition}
\smallskip

A trace arc of a generic graph may consist of several edges in the usual sense.
\smallskip
 
\begin{lemma}
\label{lem:GenericSurfaceGraph}
{\bf (a)}
For any generic surface $S$, the trace graph $\tg(S)$
 is generic in the sense of Definition~\ref{def:GenericTraceGraph}.
So the trace graph $\tg(K)$ of a generic link $K$
 is generic.
\end{lemma}
\begin{proof}
Conditions~(i)-(v) of Definition~\ref{def:GenericSurface} imply
 Conditions~(i)-(ii) of Def.~\ref{def:GenericTraceGraph}.
\end{proof}
\smallskip

\begin{figure}
\includegraphics[scale=1.0]{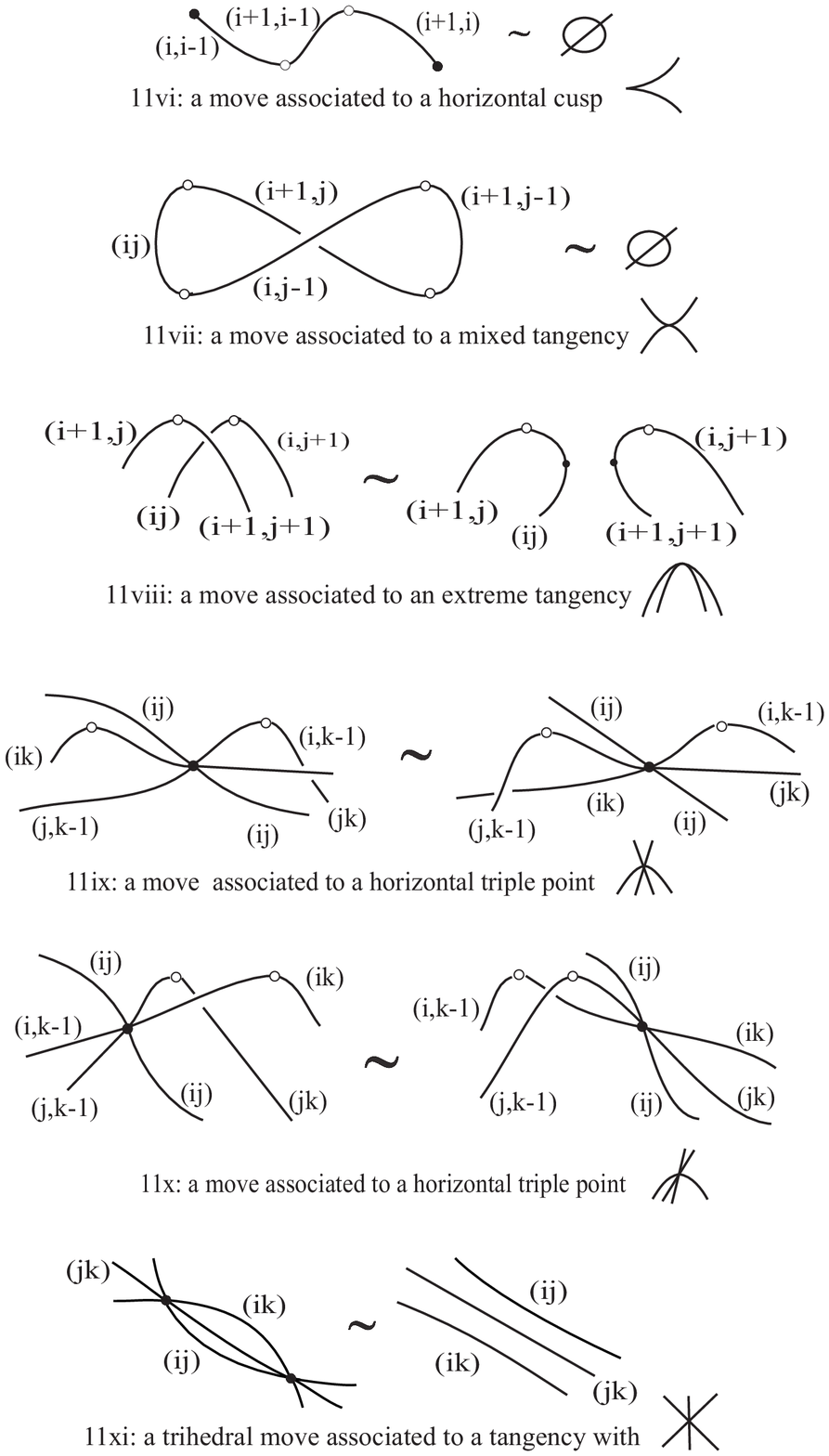}
\caption{Moves on trace graphs}
\label{fig:MovesTraceGraphs}
\end{figure}

\begin{definition}
\label{def:EquivalenceTraceGraphs}
A smooth family of trace graphs $\{G_s\}$, $s\in[0,1]$,
 is called an \emph{equivalence} if
 there are finitely many critical moments
 $s_1,\dots,s_k\in[0,1]$ such that
\smallskip

$\bu$
 for all non-critical moments $s\notin\{s_1,\dots,s_k\}$,
 the trace graphs $G_s$ are generic;
\smallskip

$\bu$
 if $s$ passes through a critical moment,
 $G_s$ changes by a move in Figure~\ref{fig:MovesTraceGraphs}.
\ed
\end{definition}
\smallskip

\begin{figure}
\includegraphics[scale=1.0]{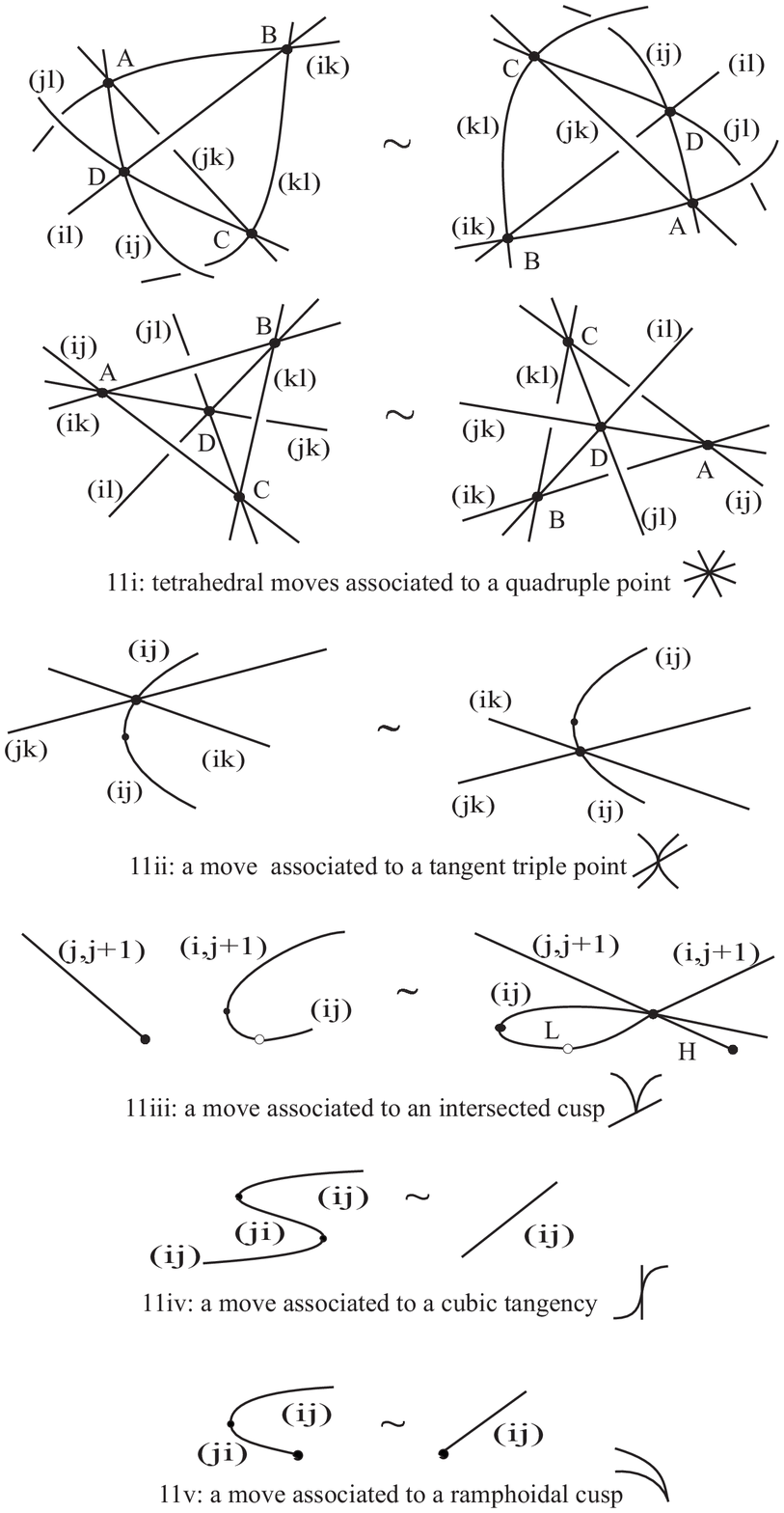}
\end{figure}

\begin{figure}[!h]
\includegraphics[scale=1.0]{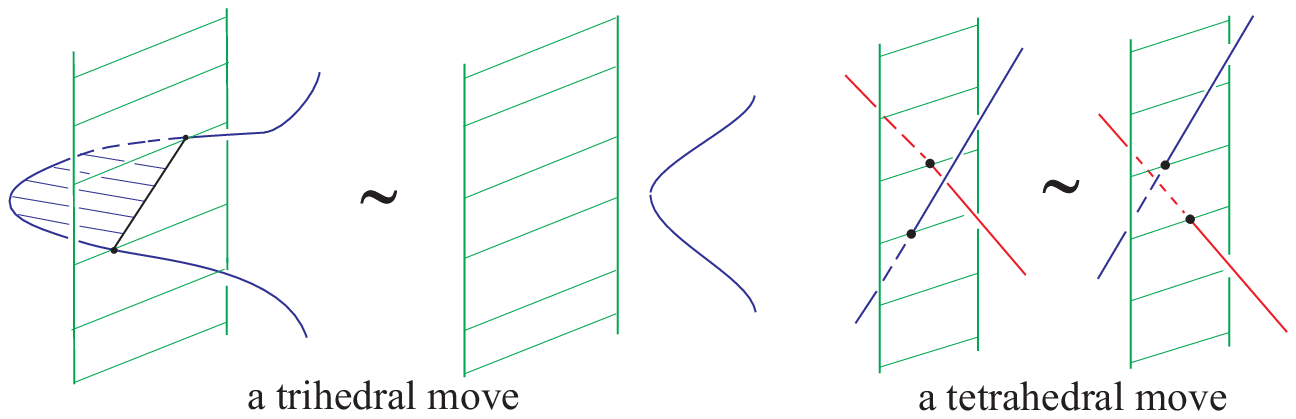}
\caption{A trihedral move and a tetrahedral move for links.}
\label{fig:MovesOnBraids}
\end{figure}

The moves in Figure~\ref{fig:MovesTraceGraphs} should be considered locally,
 ie the diagrams do not change outside the pictures.
Various mirror images of the moves are also possible.
Moreover, some labels $s+1$ can be replaced by $s-1$ and
 vice versa.
Trace graphs are symmetric under $t\mapsto t+\pi$, i.e. each move 
 in Figure~\ref{fig:MovesTraceGraphs} denotes two symmetric moves.
The most non-trivial moves are \emph{tetrahedral} moves~\ref{fig:MovesTraceGraphs}i
 and \emph{trihedral} moves~\ref{fig:MovesTraceGraphs}xi.
Their geometric interpretation at the level of links
 is shown in Figure~\ref{fig:MovesOnBraids}.

Notice that both moves in Figure~\ref{fig:MovesTraceGraphs}i 
 can be realized for links and closed braids.
In general a tetrahedral move corresponds to
 a link or a braid with a horizontal quadrisecant.
Geometrically two arcs intersect a wide band bounded by
 another two arcs.
Under a tetrahedral move, the two intersection points
 swap their heights as in Figure~\ref{fig:MovesOnBraids}.
The first picture of Figure~\ref{fig:MovesTraceGraphs}i applies when
 the intermediate oriented arcs go together
 from one side of the band to another like $\rightrightarrows$.
The second picture means that the arcs are antiparallel
 as in the British rail mark $\rightleftarrows$.
It is easier to understand Lemma~\ref{lem:ReconstructSurfaceGraph} first for knots, 
 when the indices $i,j=1$ can be missed. 

\begin{lemma}
\label{lem:ReconstructSurfaceGraph}
{\bf (a)}
For a generic trace graph $G$ such that $G\cap(\axz\times\{0\})$
 are crossings of a general diagram,
 there is a generic surface $S$ such that $\tg(S)=G$.
\medskip

\noindent
{\bf (b)}
For any equivalence of trace graphs $\{G_r\}$, there is
 an equivalence of surfaces $S_r$ with $\tg(S_r)=G_r$,
 $r\in[0,1]$.
\end{lemma}
\begin{proof}
{\bf (a)}
Consider a vertical section $P_t=G\cap(\axz\times\{t\})$
 not containing vertices of $G$.
Then $P_t$ is a finite set of points with labels
 $(q_i s_j)$, where $i,j\in\{1,\dots,m\}$, see Definition~5.2.
The points in $P_t$ will play the role of
 crossings of sections of $S$.
\smallskip

The labelled set $P_t$ defines the Gauss diagram $\gd_t$
 as follows, see Definition~\ref{def:GaussDiagram}.
Take $\sqcup_{i=1}^m S_i^1$, split each circle $S_i^1$ into $n_i$ arcs 
 and number them by $1,\dots,n_i$ according to the orientation.
We mark several points in the $q$-th arc of $S_i^1$ 
 in a 1-1 correspondence and the same order with the points of $P_t$ projected
 under $\pr_z:P_t\to S_z^1$ and having labels $(q_i s_j)$ or $(s_j q_i)$, $s=1,\dots,n_j$.
\smallskip

So each point of $P_t$ gives 2 marked points in $\sqcup_{i=1}^m S_i^1$,
 labelled with $(q_i s_j)$ and $(s_j q_i)$.
Connect them by a chord and get the Gauss diagram $\gd_t$.
The zero Gauss diagram $\gd_0$ is realizable by the given general diagram.
Hence all Gauss diagrams $\gd_t$ give rise to
 a family of diagrams $D_t$, ie to a surface $S=\cup(D_t\times\{t\})$.
\medskip

\noindent
{\bf (b)}
Apply the construction from {\bf (a)}
 to each trace graph $G_r$, $r\in[0,1]$.
\end{proof}

\begin{proposition}
\label{prop:ReductionToTraceGraphs}
{\bf (a)} 
Trace graphs $\tg(S_0),\tg(S_1)$ of generic surfaces
 are equivalent in the sense of Definition~\ref{def:EquivalenceTraceGraphs} 
 if and only if the surfaces $S_0,S_1$ are equivalent 
 in the sense of Definition~\ref{def:EquivalenceSurfaces}.
\smallskip

\noindent
{\bf (b)}
Generic surfaces $S_0,S_1$ are equivalent
 in the sense of Definition~\ref{def:EquivalenceSurfaces}
 if and only if $\tg(S_0),\tg(S_1)$ are equivalent
 in the sense of Definition~\ref{def:EquivalenceTraceGraphs}.
\end{proposition}
\begin{proof}
{\bf (a), (b)}
Any equivalence $\{S_r\}$ of surfaces gives rise to
 the equivalence $\tg(S_r)$ of trace graphs.
Any equivalence of trace graphs gives
 rise to a smooth family of diagram surfaces $\{S_r\}$
 by Lemma~\ref{lem:ReconstructSurfaceGraph}b.
The family $\{S_r\}$ is an equivalence of diagram
 surfaces since the moves in Figure~\ref{fig:MovesTraceGraphs}  are restrictions of
 the moves in Figure~\ref{fig:MovesDiagramSurfaces}.
\end{proof}
\smallskip

\noindent
{\bf Theorem~\ref{thm:MovesTraceGraphs}} directly follows from
 Propositions~\ref{prop:ReductionToGenericLinks}, \ref{prop:ReductionToGenericLoops}, 
 \ref{prop:ReductionToDiagramSurfaces} and \ref{prop:ReductionToTraceGraphs}.
\medskip

\begin{lemma}
\label{lem:ReconstructLinkGraph}
Suppose that the trace graph $G=\tg(K)$ of 
 a generic link $K$ is given, but $K$ is unknown.
Then one can construct a generic link $K'$ equivalent to $K$.
\end{lemma}
\begin{proof}
Lemma~\ref{lem:ReconstructSurfaceGraph}a 
 provides a generic surface $S$ such that $\tg(S)=G$.
Due to labels of trace arcs, the section $D_0=S\cap(\axz\times\{0\})$
 gives rise to a link $K'\subset V$ with $\pr_{xz}(K')=D_0$.
The link $K'$ can be assumed to be generic by 
 Proposition~\ref{prop:ReductionToGenericLinks}a
 and is equivalent to $K$ since
 $K$ and $K'$ have the same Gauss diagram.
\end{proof}
\smallskip


\subsection{Combinatorial construction of a trace graph}
\label{subs:Constructions}
\noindent
\smallskip

\begin{lemma}
\label{lem:Construction}
Let $K\subset V$ be a link with $2e$ extrema of the projection $\pr_z:K\to S_z^1$
 and $l$ crossings in the diagram $\pr_{xz}(K)$.
Let the extrema and intersection points from $K\cap(\dxy\times\{z=\pm 1\})$
 divide $K$ into $n$ arcs monotonic with respect to $\pr_z$.
Then $K$ is isotopic in $V$ to a link $K'$
 such that $\tg(K')$ contains $2l(n-2)$ triple vertices,
 $4(n-e-1)e$ critical vertices and $2e$ hanging vertices.
\end{lemma}
\begin{proof}
Take a generic link $K'$ smoothly equivalent to $K$ and 
 having an isotopic plane diagram.
We split $K'$ by horizontal planes into
 several horizontal slices such that each slice contains exactly
 one crossing or one extremum with respect to $\pr_z:K'\to S_z^1$.
We may assume that all maxima are above all minima,
 otherwise deform $K'$ accordingly.
To each slice we associate the corresponding
 elementary trace graph and glue them together, see 
 examples in  Figure~\ref{fig:ElementaryTraceGraphs} and
 Figure~\ref{fig:TraceGraphsExtrema}.
\smallskip

\begin{figure}[!h]
\includegraphics[scale=1.0]{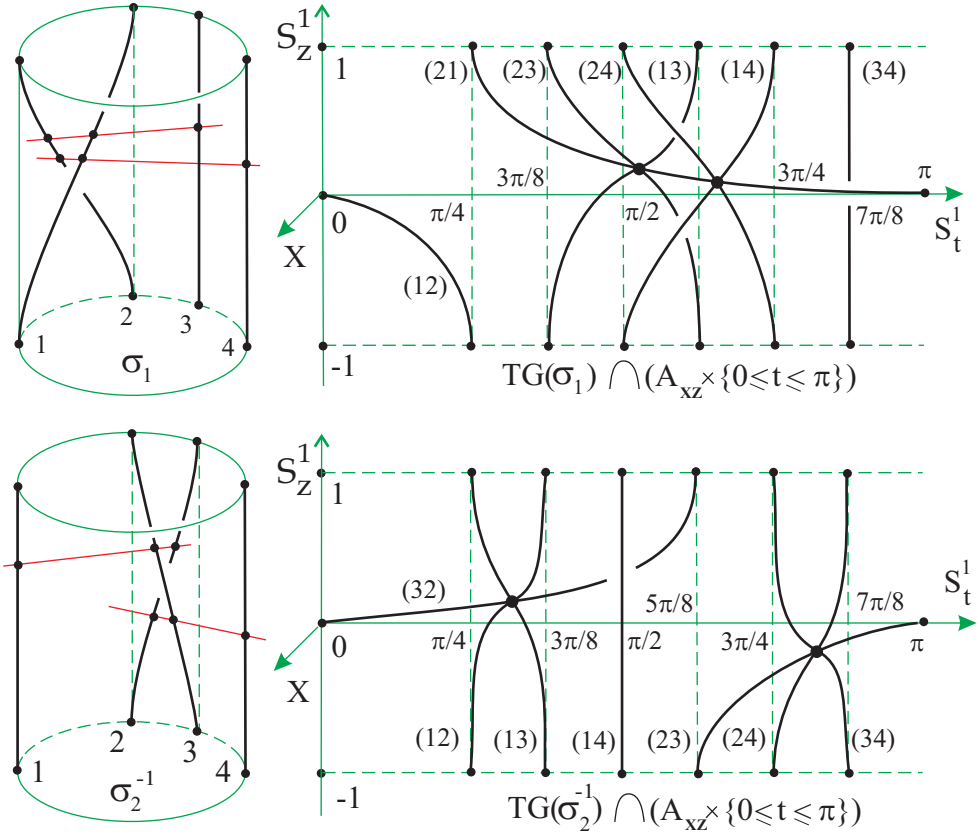}
\caption{Half trace graphs of the 4-braids
$\si_1,\si_2^{-1}\in B_4$.}
\label{fig:ElementaryTraceGraphs}
\end{figure}

\begin{figure}[!h]
\includegraphics[scale=1.0]{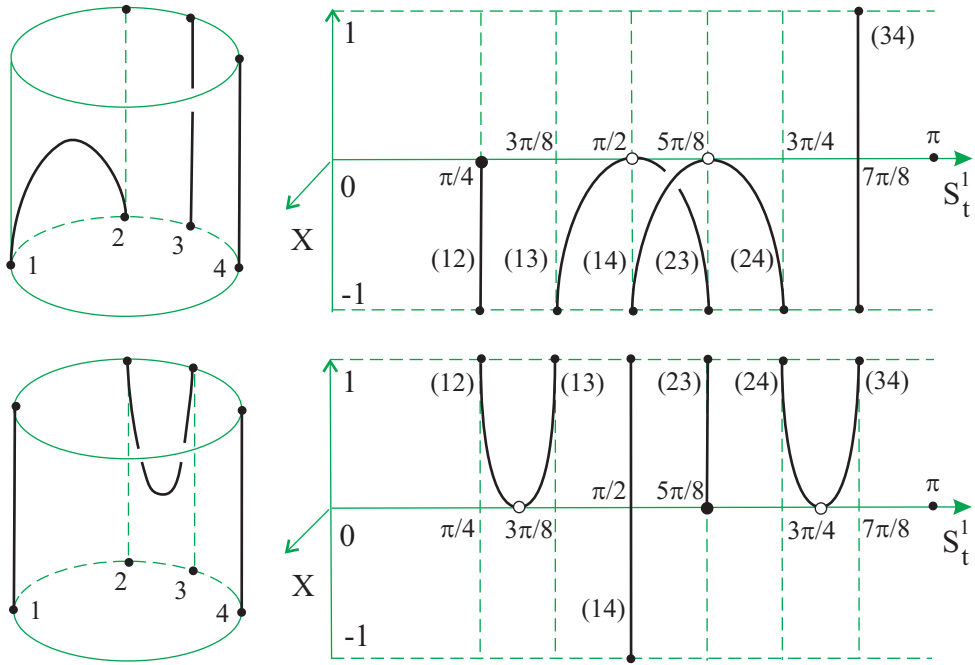}
\caption{Half trace graphs of elementary blocks containing extrema.}
\label{fig:TraceGraphsExtrema}
\end{figure}

Figure~\ref{fig:ElementaryTraceGraphs} shows two explicit examples for
 the opposite crossings in the braid group $B_4$.
In general we mark out the points $\psi_k=2^{1-k}\pi$,
 $k=0,\dots,n-1$ on the boundary of the bases $\dxy\times\{\pm 1\}$.
The $0$-th point $\psi_0=2\pi$ is the $n$-th point.
\smallskip

The crucial feature of the distribution $\{\psi_k\}$
 is that all straight lines passing through two points
 $\psi_j,\psi_k$ are not parallel to each other.
Firstly we draw all strands in the cylinder
 $\bd\dxy\times[-1,1]_z$.
Secondly we approximate with the first derivative
 the strands forming a crossing by smooth arcs,
 see the left pictures in Figure~\ref{fig:ElementaryTraceGraphs}.
\smallskip

Then each elementary braid $\si_i$ constructed as above has exactly
 $n-2$ horizontal trisecants through the strands $i,i+1$ and $j$
 for $j\neq i,i+1$.
Each trisecant is associated to a triple vertex of the trace graph,
 see 4 horizontal trisecants in the left picture of Figure~\ref{fig:ElementaryTraceGraphs}.
The trace graphs in Figure~\ref{fig:ElementaryTraceGraphs} are not generic in the sense
 of Definition~\ref{def:GenericTraceGraph}, eg
 parallel strands 3 and 4 lead to the vertical trace arc labelled with $(34)$.
But we may slightly deform such a trace graph to
 make it generic.
\smallskip

In the first picture of Figure~\ref{fig:TraceGraphsExtrema} the arc with a maximum
 is the intersection of the cylinder $\bd\dxy\times[-1,1]_z$
 with an inclined plane containing the straight line 1-2 in
 the base $\dxy\times\{-1\}$.
The highest maximum of $K'$ leads to exactly $2(n-2e)$
 critical vertices (with symmetric images under $t\mapsto t+\pi$), 
 the next maximum gives $2(n-2e+2)$ critical vertices and so on, 
 i.e. the total number is $2(n-2e)+2(n-2e+2)+\cdots+2(n-2)=2(n-e-1)e$.
The number of critical vertices associated to minima of $K'$ is the same.
Moreover each of $2e$ extrema gives one hanging vertex.
\end{proof}
\smallskip


\end{document}